\documentclass{amsart}
\usepackage{amssymb}
\usepackage{graphicx}
\usepackage{latexsym,amscd}
\usepackage{tikz}
\usepackage{tikz-cd}
\title[Some  Imaginary Quadratic Fields]
      {The Computation of  $\mathbf{\Gal(k^\infty/k)}$ for some
        Complex Quadratic Number Fields $\mathbf{k}$}

\author{E.\,Benjamin, F.\,Lemmermeyer, C.\,Snyder}


\begin{document}
\newcommand{\rsp}{\raisebox{0em}[2.7ex][1.3ex]{\rule{0em}{2ex} }}
\newcommand{\ul}{\underline}
\newcommand{\1}{\boldsymbol{1}}
\newcommand{\I}{\boldsymbol{i}}
\newcommand{\J}{\boldsymbol{j}}
\newcommand{\K}{\boldsymbol{k}}
\newcommand{\C}{{\mathbb C}}
\newcommand{\h}{\mathbb{H}}
\newcommand{\Q}{{\mathbb Q}}
\newcommand{\R}{{\mathbb R}}
\newcommand{\B}{{\mathbb B}}
\newcommand{\F}{{\mathbb F}}
\newcommand{\N}{{\mathbb N}}
\newcommand{\Z}{{\mathbb Z}}
\newcommand{\aQ}{{\overline{\Q}}}
\newcommand{\cA}{\mathcal A}
\newcommand{\OO}{\mathcal O}
\newcommand{\Rl}{{\rm Re}}
\newcommand{\cM}{\mathcal M}
\newcommand{\Am}{\operatorname{Am}}
\newcommand{\eq}{\stackrel{2}{=}}
\newcommand{\fZ}{\mathfrak 2}
\newcommand{\fQ}{\mathfrak Q}
\newcommand{\fa}{\mathfrak a}
\newcommand{\fb}{\mathfrak b}
\newcommand{\fc}{\mathfrak c}
\newcommand{\fC}{\mathfrak C}
\newcommand{\fE}{\mathfrak E}
\newcommand{\fp}{\mathfrak p}
\newcommand{\fP}{\mathfrak P}
\newcommand{\fA}{\mathfrak A}
\newcommand{\fB}{\mathfrak B}
\newcommand{\fG}{\mathfrak G}
\newcommand{\fH}{\mathfrak H}
\newcommand{\fJ}{\mathfrak J}
\newcommand{\fO}{\mathfrak O}
\newcommand{\fr}{\mathfrak r}
\newcommand{\frb}{\widetilde{\fr}}
\newcommand{\fR}{\mathfrak R}
\newcommand{\fq}{\mathfrak q}
\newcommand{\gac}{\widehat}
\newcommand{\Ft}{\widetilde{F}}
\newcommand{\disc}{\operatorname{disc}}
\newcommand{\ch}{\operatorname{char}}
\newcommand{\id}{\operatorname{id}}
\newcommand{\rank}{\operatorname{rank}}
\newcommand{\lcm}{{\operatorname{lcm}}}
\newcommand{\gen}{{\operatorname{gen}}}
\newcommand{\cen}{{\operatorname{cen}}}
\newcommand{\Tr}{\operatorname{Tr}}
\newcommand{\Cl}{\operatorname{Cl}}
\newcommand{\Ram}{\operatorname{Ram}}
\newcommand{\Gal}{\operatorname{Gal}}
\newcommand{\Aut}{\operatorname{Aut}}
\newcommand{\Inn}{\operatorname{Inn}}
\newcommand{\conj}{\operatorname{conj}}
\newcommand{\AGL}{\operatorname{AGL}}
\newcommand{\Ind}{\operatorname{Ind}}
\newcommand{\GL}{\operatorname{GL}}
\newcommand{\im}{\operatorname{im}\,}
\newcommand{\eps}{\varepsilon}
\newcommand{\la}{\langle}
\newcommand{\ra}{\rangle}
\newcommand{\lra}{\rightarrow}
\newcommand{\too}{\mapsto}
\newcommand{\ts}{\textstyle}
\newcommand{\ab}{\operatorname{ab}}
\newcommand{\ff}{\mathfrak f}
\newcommand{\fm}{\mathfrak m}
\newcommand{\fD}{\mathfrak D}
\newcommand{\fN}{\mathfrak N}
\newcommand{\fn}{\mathfrak n}
\def\la{\langle}
\def\ra{\rangle}
\def\lra{\rightarrow}
\def\lms{\mapsto}
\def\vp{\varphi}
\def\vt{\vartheta}
\def\tr{\mbox{{\rm Tr}}\,}
\def\hom{\mbox{{\rm Hom}}\,}
\def\ind{\mbox{{\rm Ind}}}
\newtheorem{df}{Definition}
\newtheorem{thm}{Theorem}
\newtheorem{lem}{Lemma}
\newtheorem{prop}{Proposition}
\newtheorem{cor}{Corollary}

\begin{abstract}
  We determine $\Gal(k^\infty/k)$ up to isomorphism for two particular
  families of imaginary quadratic number fields $k$ with $2$-class
  field tower of length $2$, where $k^\infty$ denotes the union of the
  $2$-class field tower of $k$.
\end{abstract}

\maketitle

\section{{\bf Introduction}}
Let $k$ be an algebraic number field and $\{k^n\}_{n=0}^\infty$ the
$2$-class field tower of $k$, i.e., the sequence of fields, $k^n$,
where $k^0=k$, $k^1$ is the Hilbert $2$-class field of $k$ (the
maximal everywhere-unramified abelian extension of $2$-power degree
over $k$), and $k^{n+1}=(k^n)^1$. The length of the tower is $1$ less
than the cardinality of the set $\{k^n:n\geq 0\}$. As usual, we let
$\Cl(k)$ be the ideal class group of $k$, with order $h(k)$, and
$\Cl_2(k)$ its $2$-class group of order $h_2(k)$.

In a recent article, \cite{BS}, we completed the project of
determining all imaginary quadratic fields $k$ for which $\Cl_2(k^1)$
has $2$-rank at most $2$. It follows that the $2$-class field tower of
each of these $k$ terminates at $k^j$, for $j\leq 2$.  (This fact is
trivial for $2$-$\rank\Cl_2(k^1)\leq 1$, and is described in \cite{BS}
for $2$-rank $=2$.)  This suggests that it is perhaps ``relatively
easy" to determine the Galois groups $\Gal(k^\infty/k)$ for all the
imaginary quadratic fields $k$ for which $2$-$\rank\Cl_2(k^1)\leq 2.$
This has already been done for all these $k$ with cyclic (including
trivial) $\Cl_2(k^1)$; cf.\,Table II in \cite{BLS0}.  Hence we are
left to consider $\Cl_2(k^1)$ of $2$-rank $2$. Some progress in this
case has already been made. For instance, in \cite{BLS1} we determined
$\Gal(k^\infty/k)$ for precisely those fields $k$ for which
$\Cl_2(k)\simeq (2,2,2)$ and $4$-$\rank\Cl_2(k^1)\leq 1$. In light of
\cite{BS}, this leaves open the following (not mutually exclusive)
categories for which the $2$-rank of $\Cl_2(k^1)$ is possibly $2$:

\bigskip

(I)\;\;$\Cl_2(k)\simeq (2^m,2^n),\; m,n\geq 1,$ not both $=1$\;
(cf.\,\cite{BLS2} and \cite{BS2}),

\medskip

(II)\;\; $\Cl_2(k)\simeq (2,2,2^n),\; n\geq 2$\; (cf.\,\cite{BS}),

\medskip

(III)\;\;$\Cl_2(k^1)\simeq (2^m,2^n), \;m,n\geq 2.$

\bigskip

In this present work we consider a subcase of category II: a
particular family of those $k$ for which $\Cl_2(k)\simeq (2,2,2^n)$
($n\geq 2$) such that $\Cl_2(k^1)\simeq (2,2^m)$ with $m\geq
1$. Namely, we compute $\Gal(k^\infty/k)$ for those fields $k$ of
Type~(4p) as given in Theorem~6 of \cite{BS}, i.e., those $k$ with
discriminant $d_k=-4qq'p$ where $q,q', p$ are primes such that
$q,q'\equiv 3\bmod 8$, $p\equiv 1\bmod 8$, and
$$(p/q)=(p/q')=(-q/q')=-1.$$

We then consider an obvious natural extension of the groups that occur
as the Galois groups of the fields of Type~(4p). This ends up giving
us the determination of $\Gal(k^\infty/k)$ for all the fields of
Type~(4r) (as given in Theorem~7 of \cite{BS}). See Section~\ref{S7}
below for all the details.

\bigskip

\section{\bf Structure of the Galois Groups}

\begin{thm}\label{T1}
  Let $k$ be an imaginary quadratic number field with discriminant
  $d_k=-4pqq'$ where $p,q,q'$ are primes such that
\begin{equation}\label{ET}
 p\equiv 1\bmod 8,\;q,q'\equiv 3\bmod 8,\; (p/q)=(p/q')=(-q/q')=-1.
\end{equation}
Then $\Cl_2(k)\simeq (2,2,2^n)$ and $\Cl_2(k^1)\simeq (2,2^m),$ for
some $n,m\geq 2$.

Moreover, $\Gal(k^\infty/k)\simeq \Gamma_{n,m},$ where $\Gamma_{n,m}$
is given by
$$\Gamma_{n,m}=\la a_1,a_2,a_3:
a_1^2=c_{13}^{-1},\;a_2^2=c_{13}^{2^{m-1}},\;
a_3^{2^n}=c_{12}c_{13}^{2^{m-1}},\;c_{23}=c_{12}^2=c_{13}^{2^m}=1\ra,$$
with  $2^n=h_2(k)/4$ and  $2^m=h_2(-p)$.

Finally, any imaginary quadratic field $k$ for which
$\Gal(k^\infty/k)=\Gamma_{n,m}$, for some integers $n,m>1$, has
discriminant $d_k=-4pqq'$ satisfying {\rm(\ref{ET})}.
\end{thm}

(Here, we have denoted by $(2,2,2^n)$ any group which is isomorphic to
the direct sum of two groups of order $2$ and a cyclic group of order
$2^n$, and by $c_{ij}$ the commutator $[a_i,a_j]$; cf.\;the appendix
below.)

\bigskip

\section{\bf Arithmetic of some Intermediate Fields}

The proof of this theorem will follow from the series of results
below. We start by determining some arithmetic in the seven unramified
quadratic extensions of $k$ and also by computing the $2$-class number
of the genus field $k_\gen$ of $k$. We will then use this information
to compute the Galois group of the $2$-class field tower of $k$ and verify that
it is indeed $\Gamma_{n,m}.$

We first determine the quadratic number fields $\Q(\sqrt{d}\,)$
contained in $k^1$ that have even class number and compute their
$2$-class numbers:

$$\begin{array}{c|cccccccc}
     d & -p & -pq & -pq'  & -qq' & pq & pq' & pqq' & -pqq' \\ \hline
 h_2(d) & 2^\mu & 2 & 2 & 4& 2& 2& 2& 2^{n+2}
 \end{array}\,,$$

\noindent where by genus theory, $\mu\geq 2$\;(cf.\,\cite{RR}). (Here,
$h_2(d)=h_2(\Q(\sqrt{d}\,).$)

\bigskip

We now determine the seven quadratic extensions, $k_j$, of $k$
contained in $k^1$ and compute various data concerning them for later
use:

\begin{table}[ht!]
 $$ \begin{array}{c|c|c|c|c|c|c}
  j  & k_j  & E_{k_j}  & N_{k_j/k}E_{k_j}  & \#\kappa_j & h_2(k_j)  & \kappa_j \\ \hline
  1  & k(\sqrt{-p}\,)& \la -1, \varepsilon_{qq'}\ra & 1 & 4 & 2^{n+1+\mu} & \la [\fp],[\fq]\ra \\
 2  & k(\sqrt{p}\,)& \la -1, \varepsilon_{p}\ra & \la -1\ra  & 2 & 2^{n+3} & \la [\fp]\ra \\
 3  & k(\sqrt{-1}\,)& \la i, \varepsilon_{pqq'}\ra & 1 & 4 & 2^{n+2} & \la [\mathtt{2}],[\fp]\ra \\
 4  & k(\sqrt{-q}\,)& \la \zeta_q, \varepsilon_{pq'}\ra & 1 & 4 & 2^{n+2} & \la [\mathtt{2}],[\fq]\ra \\
 5  & k(\sqrt{-q'}\,)& \la \zeta_{q'}, \varepsilon_{pq}\ra & 1 & 4 & 2^{n+2} & \la [\mathtt{2}],[\fp\fq]\ra \\
 6  & k(\sqrt{q}\,)& \la -1, \varepsilon_{q}\ra & 1 & 4 & 2^{n+2} & \la [\mathtt{2}],[\fq]\ra \\
 7  & k(\sqrt{q'}\,)& \la -1, \varepsilon_{q'}\ra & 1 & 4 & 2^{n+2} & \la [\mathtt{2}],[\fp\fq]\ra
             \end{array} $$
\caption{}
 \end{table}

Here, $E_K$ denotes the unit group in an algebraic number field $K$,
$\varepsilon_d$ is the fundamental unit in $\Q(\sqrt{d}\,)$,
$N_{k_j/k}$ the relative norm from $k_j$ to $k$, $\kappa_j$ the
subgroup of ideal classes in $\Cl(k)$ which capitulate in $k_j$, and
where $h_2(-p)=2^\mu$ ($\mu\geq 2$) and $\zeta_t$ is a generator of
the unit group of $\Q(\sqrt{-t}\,)$.  Finally, $\mathtt{2},\fp,$ and
$\fq$ are the prime ideals of $k$ containing $2,p,$ and $q$,
respectively.

\bigskip
For example, consider $k_1$. Then $\Gal(k_1/\Q)\simeq (2,2)$ and the
three quadratic subfields in $k_1$ are $F_0=k,\;F_1=\Q(\sqrt{-p}\,),$
and $F_2=\Q(\sqrt{qq'}\,)$. Hence by Kuroda's class number formula,
cf.\,\cite{Lem1}, formula (1.5), we have
$$ h_2(k_1)=\frac1{2}\,q(k_1/\Q)\,h_2(F_0)\,h_2(F_1)\,h_2(F_2), $$
where the unit index $q(k_1/\Q)=(E_{k_1}:E_{F_0}E_{F_1}E_{F_2}).$ But by our
list of quadratic fields above, $h_2(F_0)=2^{n+2},\,
h_2(F_1)=2^{\mu},$ and $ h_2(F_2)=1.$ Moreover the unit index
$q(k_1/\Q)=1$, since $k$ is complex quadratic and $k_1/k$ is
unramified; see, e.g.,\,\cite{Lem3}, Theorem~1. This also shows that
$E_{k_1}=\la -1,\varepsilon_{qq'}\ra$, which then implies that
$N_{k_1/k}E_{k_1}=\la 1\ra$.  Thus we have $h_2(k_1)=2^{n+1+\mu}$.
But then \\ $\#\kappa_1=2(E_k:N_{k_1/k}E_{k_1})=4$, cf.\,\cite{Iwa}
and Chapitre IV in \cite{Che} (or \cite{Ros}). Finally, since the
prime ideal above $p$ in $\Q(\sqrt{-p}\,)$ is principal, we see that
$\fp\mathcal{O}_{k_1}$ is also principal. Similarly, since
$h_2(qq')=1$, the prime ideal in $\Q(\sqrt{qq'}\,)$ above $q$ is
principal, implying that $\fq\mathcal{O}_{k_1}$ is principal,
too. Hence $[\fp]$ and $[\fq]$ are in $\kappa_1$. Since
$\#\kappa_1=4$, we thus have $\kappa_1=\la [\fp],[\fq]\ra.$ The rest
(with the exception of showing that $[\mathtt{2}]\in\kappa_j$ for
$j=4,5$ and $[\fp]\in \kappa_3$, which follows from the lemma below)
are computed in a similar manner. (It should be observed that in a
number of places we have used the fact that
$[\fp][\fq][\fq']=[(\sqrt{pqq'}\,)]=1$, where $\fq'$ is the prime
ideal of $k$ containing $q'$.) See \cite{Lem}, pp.\,206-208, for more
details.

\begin{lem}\label{L1}
  (i) Let $F_1$ be a quadratic number field where the discriminant
  $d_{F_1}=4qp$ with $q\equiv 3\bmod 8$, $p\equiv 1\bmod 8$, and
  $(p/q)=-1.$ Then $h_2(F_1)=2$ and the prime ideal $\mathtt{2}$
  containing $2$ is principal.

     (ii) Let $F_2$ be a quadratic number field such that the
  discriminant $d_{F_2}=pqq'$ with primes $p,q,q'$ satisfying $p\equiv
  1\bmod 4$, $q,q'\equiv 3\bmod 4$, and $(p/q)=(p/q')=-1$. Then
  $h_2(F_2)=2$ and the prime ideal $\fp$ containing $p$ is principal.
\end{lem}
\begin{proof}
  (i) By genus theory, $h_2(F_1)=2.$ Next, consider the genus
  characters \\ $\chi_{-4}, \chi_{-q},\chi_{p}$ and notice that
  $\chi(\mathtt{2}\fq\fp)=+1 $ for each of them (and where $\fq,\fp$
  are the prime ideals of $F_1$ containing $q,p$, respectively). This
  implies that the narrow class $[\mathtt{2}\fq\fp]_+$ is the square
  of a narrow class in $\Cl_2(F_1)^+$. (See \cite{Lem2}, subsections
  2.1-2.3\,.)  Hence, since $\Cl_2(F_1)^+\simeq (2,2)$ (again by genus
  theory), we see $\mathtt{2}\fq\fp$ is principal in the narrow sense
  and thus principal in the ordinary sense, too. But notice that
  $\fq\fp$ is principal. Thus so is $\mathtt{2}.$

The proof of (ii) is similar and left to the reader.
\end{proof}

Next, consider our field $k$ again. Then $\mathtt{2},\fp,\fq,\fq'$ are
the prime ideals of $k$ that ramify. It thus follows that the subgroup
of ideal classes of orders 1 or 2, $\Cl_2(k)[2]$, is equal to $\la
[\mathtt{2}],[\fp],[\fq],[\fq']\ra=\la [\mathtt{2}],[\fp],[\fq]\ra,$
since $\fp\fq\fq'\sim 1$. Now, since $\Cl_2(k)\simeq (2,2,2^n)$, with
$n\geq 2$, there is exactly one nontrivial class in $\Cl_2(k)[2]$ that
is a square of an ideal class.  To determine it, consider the four
genus characters associated with the prime discriminants, $d_j$, of
$d_k$: $\chi_{-4},\chi_{-q},\chi_{-q'}, \chi_p$. Then it can be seen
that $\chi_{d_j}(\mathtt{2}\fp)=+1$ for each $j$. Hence
$[\mathtt{2}\fp]$ is the unique ideal class of order $2$ which is a
square of an ideal class. Equivalently,
\begin{equation}\label{E1}
    \Cl_2(k)^2\cap \Cl_2(k)[2]=\la [\mathtt{2}\fp]\ra.
\end{equation}

 We will use this information later, when determining a presentation
 of $\Gal(k^2/k)$.

Now we come to what turns out to be the main obstacle in determining
the structure of $G=\Gal(k^2/k)$. We know at this point that $G'\simeq
(2,2^m)$ for some integer $m\geq 1$ (we'll soon see $m\geq 2$), but we
do not know if $h_2(-p)=2^m$, i.e.,\,whether or not $\mu=m$ where
$2^\mu=h_2(-p)$. As we will see later, this problem will be resolved
by proving that $h_2(k_\gen)=\frac1{4}h_2(k)h_2(-p)=2^{n+\mu}$. We
show this by using the fact that $k_\gen$ is a CM-field whose maximal
real subfield $k_\gen^+$ has odd class number, as follows from the
next lemma. (The idea here is similar to one in the proof of Theorem~3
in \cite{BLS1}.)

\begin{lem}\label{L2}
  Let $M=\Q(\sqrt{p},\sqrt{q},\sqrt{q'}\,)$, where $p,q,q'$ are
  distinct primes such that
  $$ p\equiv 1\bmod 4,\; q\equiv q'\equiv 3\bmod 4,\;(p/q)=(p/q')=-1. $$
Then $M$ has odd class number.
\end{lem}

\begin{proof}
  Let $K$ be the subfield of $M$ given by
  $K=\Q(\sqrt{pq},\sqrt{q'}\,).$ We claim that $h_2(K)=2.$ To see
  this, first notice that $\Q(\sqrt{pq}\,), \Q(\sqrt{q'}\,),$ and
  $\Q(\sqrt{pqq'}\,)$ are the intermediate quadratic fields with unit
  groups $E_{\Q(\sqrt{d}\,)}$, denoted by $e_d$, equaling $\la
  -1,\varepsilon_d\ra$, with fundamental unit $\varepsilon_d$.  Then
  $$h_2(K)=\frac1{4}q(K/Q)h_2(pq)h_2(q')h_2(pqq').$$
  Given the relations among $p,q,q'$ as stated in the lemma, we see by
  genus theory that $h_2(q')=1$ and $h_2(pq)=h_2(pqq')=2$.  Thus
  $h_2(K)=q(K/\Q)$.  Now, $q(K/\Q)=(E_K:e_{q'}e_{pq}e_{pqq'})=(E_K:\la
  -1,\varepsilon_{q'},\varepsilon_{pq},\varepsilon_{pqq'}\ra).$ Then
  by applying, for example, Proposition~3 in \cite{BLS}, to each of
  the fundamental units, we find
  $E_K=\la -1,\varepsilon_{pqq'},\varepsilon_{q'},
  \sqrt{\varepsilon_{pq}\varepsilon_{q'}}\,\ra$, and so $q(K/\Q)=2.$
  This implies that $h_2(K)=2 $, as claimed.

From this we see that $M=K^1$, the $2$-class field of $K$, and
therefore $h_2(M)=1$, as desired.
\end{proof}

\begin{prop}\label{P1}
  Let $k$ be as described in Theorem~\ref{T1} and let $k_\gen$ be its
  genus field. Then
  $$ h_2(k_\gen)=\frac1{4}h_2(k)h_2(-p)=2^{n+\mu}, $$
  where $2^n=h_2(k)/4$ and $2^\mu=h_2(-p).$
\end{prop}

\begin{proof}
  Let $L=k_\gen.$ We have
  $L=k_\gen=\Q(i,\sqrt{p},\sqrt{q},\sqrt{q'}\,).$ This is a complex
  composite of quadratic number fields whose degree over $\Q$ is
  $2^4=16.$ Hence by Kuroda's multiquadratic class number formula, see
  \cite{BLS1}, formula (1), for example,
  $$ h_2(L)=\frac1{2^{16}}q(L/\Q)\prod_d h_2(d), $$
  where the product is over all the quadratic number fields in $L$ and \\
  $q(L/\Q)=(E_L:\prod e_d)$, where $\prod e_d$ is the composite of
  the unit groups of all these quadratic fields. By the list of
  $h_2(d)$ above, we obtain
  $$ h_2(L)=2^{n+\mu-7}q(L/\Q).$$
  Now $L$ is a CM-field with maximal real subfield
  $M=L^+=\Q(\sqrt{p},\sqrt{q},\sqrt{q'}\,)$. By Proposition~1 in
  \cite{BLS1} we have $q(L/\Q)=Q(L)q(M/\Q).$ Hence to compute
  $q(L/\Q)$, we need only determine $Q(L)$ and $q(M/\Q)$.

  The field $M$ is a composite of the real quadratic number fields of
  $L$ and \\$[M:\Q]=8$. Hence by Kuroda's formula (again see \cite{BLS1})
  $$ h_2(M)=\frac1{2^9}q(M/Q)\prod_{d>0}h_2(d)=2^{-6}q(M/\Q), $$
  analogously to the above case.  Lemma~\ref{L2} then gives us
  $h_2(M)=1$ and so $q(M/\Q)=2^6$.

  Finally, we have $L=M(i)$ and thus $L/M$ is not essentially ramified
  (see \cite{Lem3}, p.\,349). Moreover $2\mathcal{O}_M=\fb^2$ with
  $\fb$ principal. To see this, let
  $K=\Q(\sqrt{q},\sqrt{q'}\,)\subseteq M.$ Then since $2$ is ramified
  in $\Q(\sqrt{q}\,)$ and splits in $\Q(\sqrt{qq'}\,)$ (for $qq'\equiv
  1\bmod 8$), $2\mathcal{O}_K=\mathtt{2}_1^2\mathtt{2}_2^2$ for
  distinct ideals $\mathtt{2}_1, \mathtt{2}_2$. But $K$ has odd class
  number, for by Kuroda's class number formula,
  $$ h_2(K)=\frac1{4}q(K/\Q)h_2(q)h_2(q')h_2(qq')=\frac1{4}q(K/\Q), $$
  for intermediate quadratic fields
  $\Q(\sqrt{q}\,),\Q(\sqrt{q'}\,),\Q(\sqrt{qq'}\,)$ with odd class
  numbers. But since
  $\Q(\sqrt{\varepsilon_q}\,)=\Q(\sqrt{2},\sqrt{q}\,)\nsubseteq K$, we
  see as above that $q(K/\Q)\leq 4$, whence $q(K/\Q)=4$ and
  $h_2(K)=1.$ Thus, since the ideal class $[\mathtt{2}_1\mathtt{2}_2]$
  has order at most $2$, we see that $\mathtt{2}_1\mathtt{2}_2$ is
  principal, say $\mathtt{2}_1\mathtt{2}_2=\beta\mathcal{O}_K$. Whence
  $2\mathcal{O}_M=\fb^2$ with $\fb=\beta\mathcal{O}_M$, which is
  principal.  From all this we see that, by Theorem~1 of \cite{Lem3},
  $Q(L)=2$.

  Therefore, $$h_2(k_\gen)=2\cdot 2^6\cdot
  2^{n+\mu-7}=2^{n+\mu}=\frac1{4}h_2(k)h_2(-p),$$ as desired.
\end{proof}

\bigskip

From this proposition, it follows that $G'\simeq (2,2^m)$ with $m\geq
2$. Namely, we need only show that $h_2(k^1)\geq 8.$ Now, we have
$[k_\gen:k]=8$ and so $[k^1:k_\gen]=2^{n-1}$. Thus
\begin{equation}\label{E1.5}
[k^2:k^1]\geq [k_\gen^1:k^1]=\frac{[k_\gen^1:k_\gen]}{[k^1:k_\gen]}=2^{\mu+1}\geq 8,
\end{equation}
since $\mu\geq 2$. Therefore $h_2(k^1)\geq 8.$

\bigskip

\section{\bf Some Group Theory Associated with the Galois Groups}

Refer to the appendix below for relevant information and definitions
that are used in this section. We need to show that $G$ has the
presentation given in the statement of Theorem~\ref{T1}.  For starters
we know by \cite{BS} that $G$ is metabelian, and moreover we have
$G=\la a_1,a_2,a_3\ra$ where
$$a_1^2\equiv a_2^2\equiv a_3^{2^n}\equiv 1\bmod G',\;\;n\geq 2,$$
and thus notice that $\exp(G_2/G_3)=2.$ This follows because $G_2=\la
c_{12},c_{13},c_{23},G_3\ra$ (cf.\,\cite{Hup}, p.\,258), and since
$a_1^2,a_2^2\in G_2$, by Proposition~\ref{P3} in the appendix below,
$$ c_{ij}^2\equiv c_{ij}^2c_{iji}=[a_i^2,a_j]\equiv 1\bmod G_3,\;\;
   (i=1,2; j=1,2,3).$$

By Theorem~6 in \cite{BS} we can find generators $a_1,a_2,a_3$ of $G$
satisfying the congruences above such that
$$H=\la a_1^4,a_2^2,a_3^4,a_1^2c_{13},c_{12},c_{23},G_3\ra,$$
where $H=G^{(2,2)}=(G^2)^2[G,G^2]$, since $G/H$ is of type 32.012, the
12th group of order $32$ in \cite{SW}. By Proposition~3 of \cite{BS}
we see that
$$ J=\la a_1^4,a_2^2,a_3^{2^n}, a_1^2c_{13},c_{12},c_{23},
   G_3\ra=\la a_3^{2^n},G_3\ra,$$
where $J=H\cap G_2$. This implies, in particular, that

\bigskip

\begin{equation}\label{E2}
a_1^2c_{13}\equiv a_2^2\equiv a_3^{2^n}\equiv 1\bmod J.
\end{equation}

\bigskip

The following lemma will help simplify some of the argument.

\begin{lem}\label{L3} $G_3=\la c_{13}^2\ra.$
\end{lem}
\begin{proof}  First,   $G_2=\la c_{12},c_{13},c_{23},G_3\ra$ and so
$$G_3=\la c_{12i},c_{13j},c_{23\ell},G_4: \,i,j,\ell=1,2,3\ra.$$
But now notice that $[J,G]=G_4$. This can be seen as follows:
Clearly $G_4\subseteq [J,G]$. For the reverse inclusion, one can
show by induction that $[G,G^{2^\ell}]\subseteq G_{\ell+2}$ and thus
$[G, a_3^{2^n}]\subseteq G_{n+2}\subseteq G_4$. Hence
$[J,G]\subseteq G_4$. Now since $c_{12},c_{23}\in J$, we have
$$c_{12i}\equiv c_{23j}\equiv 1\bmod G_4,\;i,j=1,2,3.$$
Therefore, $G_3$ is reduced to
$$G_3=\la c_{131},c_{132},c_{133},G_4\ra.$$
Also, since $a_1^2\equiv c_{13}\bmod J$, we see that
$c_{13j}\equiv [a_1^2,a_j]=c_{1j}^2c_{1j1}\bmod G_4$.
Hence,
$$c_{131}\equiv 1,\;c_{132}\equiv c_{12}^2,\;c_{133}\equiv c_{13}^2\bmod G_4.$$

By Proposition~\ref{P3} and the fact that $G_4\supseteq G''$
(for us $G''=1$), cf.\,\cite{Hup}, p.\,265, it follows that
$1\equiv c_{123}c_{231}c_{132}\bmod G_4$. Thus  from above  we have
$c_{12}^2\equiv c_{132}\equiv 1\bmod G_4.$ This all implies that
$$G_3=\la c_{13}^2,G_4\ra.$$
But then we have $G_4=\la c_{131}^2,c_{132}^2,c_{133}^2\ra,
G_5=\la c_{13}^4,G_5\ra$, whence
$$G_3=\la c_{13}^2,G_5\ra.$$
Continuing this process shows that $G_3=\la c_{13}^2,G_\ell\ra$ for all
$\ell\geq 3.$ But then since $G$ is nilpotent and thus $G_\ell=1$ for
sufficiently large $\ell$, we have $G_3=\la c_{13}^2\ra.$
\end{proof}

From this lemma, it follows that $G_j=\la c_{13}^{2^{j-2}}\ra$, for any $j\geq 3$.

\bigskip
Next, we look at the structure of  $G_2$.
\begin{lem}\label{L4} The generators $a_1,a_2,a_3$ of $G$ may be chosen
  so that they satisfy the properties above and such that
  $$ G_2=\la c_{13},\gamma\ra,$$ with $\gamma=c_{12}$ or $c_{23}$ and,
  moreover, so that if $\gamma=c_{12}$, then $c_{23}\in G_3$; and if
  $\gamma=c_{23}$, then $c_{12}\in G_3.$
\end{lem}
\begin{proof}

We already know that $G_2$ has rank $2$ and we have $G_3=G_2^2$, which
implies that $G_2/G_3\simeq (2,2)$. Also, since $G/H\simeq 32.012$, we
see that the order of $(G/H)_2$ is $2$.  But $(G/H)_2\simeq G_2H/H$
(cf.\,\cite{Hup}, p.\,262) and thus $(G_2H:H)=2$. But then
$$G_2/J=G_2/G_2\cap H\simeq G_2H/H,$$
and so $(G_2:J)=2.$ Therefore $(J:G_3)=2$. Also $c_{13}\not\in J$
(otherwise, there is a contradiction to the minimal presentation of
$G/H$ as given in \cite{SW}), but we know that $c_{12},c_{23}\in
J$. Hence $J=\la c_{12},c_{23},G_3\ra$, for if not, then
$c_{12},c_{23}\in G_3$, which implies that $J=G_3$; and this is not
the case. It follows that there are three possibilities:

\bigskip
(i) $c_{12}\in G_3$, $c_{23}\not\in G_3$,\;(ii) $c_{12}\not\in G_3$,
$c_{23}\in G_3$, and (iii) $c_{12}\not\in G_3$, $c_{23}\not\in G_3$.
\bigskip

Suppose (iii) is satisfied; then $c_{12}c_{23}\equiv 1\bmod G_3$. We
make the following change of generators:
$$a_1'=a_1,\;a_2'=a_2,\;a_3'=a_3a_1,$$
in which case
$$c_{12}'=[a_1',a_2']=c_{12},\;c_{13}'=c_{13},\; c_{23}'\equiv
  c_{12}c_{23}\equiv 1\bmod G_3.$$
Since this change in generators does not alter any of the previous
relations, it follows that we are reduced to (ii). This establishes
the lemma.
\end{proof}

Before summarizing what we have so far, we can add
$$c_{23}^2\equiv 1\bmod G_4;$$
for since $[J,G]=G_4$ and since $a_2^2\in J$ (and $c_{232}\equiv
1\bmod G_4$ from the proof of Lemma~\ref{L3}), we see that
$$1\equiv [a_2^2,a_3]=c_{23}^2c_{232}\equiv c_{23}^2\bmod G_4.$$

We now have

\begin{eqnarray}\label{E3}
 G_2=\la c_{13},\gamma\ra,\;G_3=\la c_{13}^2\ra,\; J=\la \gamma,c_{13}^2\ra=\la a_3^{2^n}, c_{13}^2\ra, \;(\gamma\in\{c_{12},c_{23}\}) \\
\nonumber  c_{12i}\equiv c_{23j}\equiv c_{131}\equiv c_{132}\equiv c_{23}^2\equiv c_{12}^2\equiv 1,\;c_{133}\equiv c_{13}^2\bmod G_4,\;(i,j=1,2,3).
\end{eqnarray}
 Moreover, $\gamma=c_{12}\Rightarrow c_{23}\in G_3$;\;$\gamma=c_{23}\Rightarrow c_{12}\in G_3$.

\bigskip
To obtain further information on the structure of $G'=G_2$, we study
some properties of the seven maximal subgroups $H_j$ of $G$. Here is a
table of the $H_j$ along with their commutator subgroups:

\begin{table}[ht!]
$$\begin{array}{c|l|l}
  j & H_j & H_j' \\\hline
  1 & \la a_1,a_2,a_3^2,G_2\,\ra & \la c_{12},c_{13}^{2^\lambda}\,\ra\;(\text{some}\;\lambda\geq 2) \\
  2 & \la a_2,a_3,G_2\,\ra & \la c_{23},c_{13}^2\,\ra \\
  3 & \la a_1a_3,a_2,G_2\,\ra & \la c_{12}c_{23},c_{13}^2\,\ra \\
   4 & \la a_1,a_3,G_2\,\ra & \la c_{13}\,\ra \\
  5 & \la a_1a_2,a_3,G_2\,\ra & \la c_{13}c_{23},c_{13}^2\,\ra \\
  6 & \la a_2a_3,a_1,G_2\,\ra & \la c_{12}c_{13},c_{13}^2\,\ra \\
   7 & \la a_1a_2,a_2a_3,G_2\,\ra & \la c_{12}c_{13}c_{23},c_{13}^2\,\ra.
\end{array}$$
\caption{}
\end{table}

\bigskip
For example, notice that for $H_1'$, \;$[a_1,a_2]=c_{12}$, but $[a_1,a_3^2], [a_2,a_3^2], [a_1,G_2], [a_2,G_2]$, and $[a_3^2,G_2]$ are all contained in $G_4=\la c_{13}^4\ra$, as the reader can readily verify by considering (\ref{E3}) above. The rest are computed similarly.

\bigskip

Now notice that $(H_j:G')=2^{n+1}$ and
$$(G':H_j')=\left\{
              \begin{array}{rl}
                2, & \;\;\hbox{if}\;\; j=3,4,5,6,7, \\
                2, &  \;\;\hbox{if}\;\; j=2, \;\text{and}\;\; c_{23}\not\in G_3, \\
                4, &  \;\;\hbox{if}\;\; j=2, \;\text{and}\;\; c_{23}\in G_3, \\
                \geq 2^2, & \;\;\hbox{if}\;\; j=1.
              \end{array}
            \right.$$
which implies that
$$ (H_j:H_j')=\left\{
              \begin{array}{rl}
                2^{n+2}, & \;\;\hbox{if}\;\; j=3,4,5,6,7, \\
                2^{n+2}, &  \;\;\hbox{if}\;\; j=2, \;\text{and}\;\; c_{23}\not\in G_3, \\
                2^{n+3}, &  \;\;\hbox{if}\;\; j=2, \;\text{and}\;\; c_{23}\in G_3, \\
                \geq 2^{n+3}, & \;\;\hbox{if}\;\; j=1.
              \end{array}
            \right.$$

\bigskip
Each of these $H_j$ is the subgroup fixing a unique unramified quadratic extension $k_{i_j}$ of $k$, i.e., for each $j=1,\dots,7$, there is a unique $i=i_j\in\{1,\dots,7\}$ such that $H_j=\Gal(k^2/k_{i_j}).$  We can thus compare Tables 1 and 2; in particular we have $(H_j:H_j')=h_2(k_{i_j}).$ Hence, $c_{23}\in G_3$; for if not, then six of the $H_j$ would satisfy $(H_j:H_j')=2^{n+2}$, which cannot be the case  since $h_2(k_i)=2^{n+2}$ for only five of the $k_i$, as can be seen in Table~1. Lemma~\ref{L4} then implies that $G_2=\la c_{13},c_{12}\ra$.

Also, by comparing the tables, we see that if $\mu=2$, then\\ $\{H_1,H_2\}=\{\Gal(k^2/k_1),\Gal(k^2/k_2)\},$ and if $\mu>2$, then $H_1=\Gal(k^2/k_1)$. In either case, we have $\lambda=\mu$ and then since $(G_2:H_1')=2^\mu$, $c_{12}^2\in \la c_{13}^{2^\mu}\,\ra.$ \;\; Table 2 now  becomes:

\begin{table}[ht!]
$$\begin{array}{c|l|l}
  j & H_j & H_j'\\\hline
   1 & \la a_1,a_2,a_3^2,c_{13},c_{12}\,\ra & \la c_{12},c_{13}^{2^\mu}\,\ra \\
  2 & \la a_2,a_3,c_{13},c_{12}\,\ra & \la c_{13}^2\,\ra \\
  3 & \la a_1a_3,a_2,c_{13},c_{12}\,\ra & \la c_{12},c_{13}^2\,\ra \\
   4 & \la a_1,a_3,c_{13},c_{12}\,\ra & \la c_{13}\,\ra \\
  5 & \la a_1a_2,a_3,c_{13},c_{12}\,\ra & \la c_{13}\,\ra \\
  6 & \la a_2a_3,a_1,c_{13},c_{12}\,\ra & \la c_{12}c_{13},c_{13}^2\,\ra \\
   7 & \la a_1a_2,a_2a_3,c_{13},c_{12}\,\ra & \la c_{12}c_{13},c_{13}^2\,\ra.
\end{array}$$
\caption{}
\end{table}

\bigskip
Summarizing, we now have

\begin{eqnarray}\label{E4}
 G_2=\la c_{12},c_{13}\ra,\;G_3=\la c_{13}^2\ra,\; J=\la c_{12},c_{13}^2\ra=\la a_3^{2^n}, c_{13}^2\ra, \\
\nonumber  c_{12i}\equiv c_{23j}\equiv c_{131}\equiv c_{132}\equiv 1,\;c_{133}\equiv c_{13}^2\bmod G_4,\;(i,j=1,2,3),\\
\nonumber c_{23}\equiv 1\bmod G_3,\;c_{12}^2\equiv 1\bmod \la c_{13}^{2^\mu}\ra.
\end{eqnarray}

By (\ref{E2}) and (\ref{E4}), we have
$$a_1^2c_{13}=c_{13}^{2u_1}c_{12}^{v_1},\;\;a_2^2=c_{13}^{2u_2}c_{12}^{v_2},\;\;a_3^{2^n}=c_{13}^{2u_3}c_{12}^{v_3},$$
for some integers $u_i,v_i$, with $v_3$ odd (otherwise $J\subseteq G_3$, which is not the case).

Now we can simplify the exponents here a little. For if $v_2$ is even, then $c_{12}^{v_2}\in G_4=\la c_{13}^4\ra$ and thus we may take $v_2=0$ by altering $u_2$, since $c_{12}^2\in \la c_{13}^2\ra$. Now suppose $v_2$ is odd. Then by  changing  generators to $a_1'=a_1$, $a_2'=a_2a_3^{2^{n-1}}$, and $a_3'=a_3$,  we claim we can  assume that $v_2=0$. To establish this claim, first we have
$$(i)\;\;c_{13}'=c_{13},\quad
(ii)\;\;c_{12}'=c_{12}c_{13}^{2^n\eta},\quad
(iii) \;c_{23}'=c_{23}c_{13}^{2^n\epsilon}$$
for some integers $\eta,\epsilon$. Clearly $(i)$ is obvious. With respect
to $(ii)$ we have $c_{12}'=c_{12}[a_1,a_3^{2^{n-1}}][c_{12},a_3^{2^{n-1}}]$.
But by (\ref{E4}),  $[a_1,a_3^2]=c_{13}^2c_{133}\in G_4=\la c_{13}^4\ra$
and $[a_2,a_3^2]\in G_4$, and trivially $[a_3,a_3^2]=1$. Thus by induction
\begin{eqnarray}\label{E5}
[G,a_3^{2^\ell}]\subseteq \la c_{13}^{2^{\ell+1}}\ra,\;\;\text{for any}\;\ell\geq 1.
\end{eqnarray}
We note that $(iii)$ is handled in the same way and left to the reader.

From this,

$$a_1'^2c_{13}'=a_1^2c_{13}=c_{13}^{2u_1}c_{12}^{v_1}=c_{13}'^{2u_1-2^nv_1\eta}c_{12}'^{v_1}=c_{13}^{2u_1'}c_{12}^{v_1},$$
where $u_1'=u_1-2^{n-1}v_1\eta.$ In an analogous way we get
$$a_3'^{2^n}=c_{13}'^{2u_3'}c_{12}'^{v_3}.$$

Finally,
$$a_2'^2=(a_2a_3^{2^{n-1}})^2=a_2^2a_3^{2^n}[a_3^{2^{n-1}},a_2][a_3^{2^{n-1}},a_2,a_3^{2^{n-1}}]=a_2^2a_3^{2^{n}}c_{13}^{2^n\theta},$$
for  some integer $\theta$ by (\ref{E5}). Whence,
$$a_2'^2=c_{13}'^{2u_2+2u_3-2^n(v_2+v_3)\eta+2^n\theta}c_{12}'^{v_2+v_3}.$$
But $v_2+v_3$ is even, and thus as in the previous case we may write
$a_2'^2=c_{13}'^{2u_2'}$ for an appropriate choice of $u_2'$.

We claim that the new generators  still satisfy all the previous relations in (\ref{E4}). We need only check that $c_{12}'^2\equiv c_{13}'^{2^\mu}$, as the rest are immediate. To see this, it follows  that $H_1=\la a_1',a_2',a_3'^2,c_{13}',c_{12}'\ra$,
and that thus $H_1'\supseteq \la c_{12}',c_{13}^{2^\mu}\ra.$ We have $(G_2:H_1')=2^\mu$, and thus if $c_{12}'^2\not\in\la c_{13}^{2^\mu}\ra,$ then $(G_2:H_1')\leq (G_2:\la c_{12}', c_{12}'^2\ra)< 2^\mu$, a contradiction.

Also, without loss of generality, we may assume that $v_1=0$ or $1$ and $v_3=1$, by altering $u_1$ and $u_3$, respectively.  We therefore have (again without loss of generality)
\begin{eqnarray}\label{E6}
a_1^2c_{13}=c_{13}^{2u_1}c_{12}^{v_1}\,,\;\;a_2^2=c_{13}^{2u_2}\,,\;\; a_3^{2^n}=c_{13}^{2u_3}c_{12},\;\; v_1\in\{0,1\}.
\end{eqnarray}

We now wish to obtain more precise information on the exponents $u_i$ and $v_1$. To do so, we  will need to compute transfer maps, $t_{H_i}$, from $G/G'$ to $H_i/H_i'$, then look at their kernels, and compare these with the capitulation kernels, $\kappa_j$ in Table~1. (Recall that the homomorphism $j_t:\Cl_2(k)\rightarrow \Cl_2(k_t)$ induced by extension of ideals correspond via the Artin reciprocity law to $t_{H_i}$, where $H_i=\Gal(k^2/k_t).$)

First recall that if $G$ is a finite $2$-group and $H$ is a maximal subgroup of $G$, then the transfer map, $t_H$, is the homomorphism $t_{G,H}:G/G'\rightarrow H/H'$, given as follows: Let $G=H\cup Hz$. Then
$$t_{G,H}(\overline{x})=\left\{
    \begin{array}{ll}
     x^2\,[x,z]\,H', & \hbox{if}\;\; x\in H, \\
      x^2\,H', & \hbox{if not,}
    \end{array} \right.$$
where for any $x\in G$, $\overline{x}=xG'$ in the factor group $\overline{G}=G/G'$. Here is a table of some values of the transfer maps, where we  let  $t_j=t_{G,H_j}$. Ultimately, in trying to determine $\ker t_j$, we are interested in the values in $\overline{G}[2]$, for since $(G:H_j)=2$, it follows, as is easily seen, that $\ker t_j\subseteq \overline{G}[2]$ (the subgroup of $\overline{G}$ of elements of order at most $2$). Notice that $\overline{G}[2]=\la \overline{a_1},\overline{a_2},\overline{a_3}^{2^{n-1}}\,\ra.$

\begin{table}[ht!]
$$\begin{array}{c|llll}
    j & t_j(\overline{a_1}) & t_j(\overline{a_2}) & t_j(\overline{a_3}^{2^{n-1}}) & \bmod H_j' \\\hline
    1 & c_{13}^{2u_1} & c_{13}^{2u_2}c_{23} & c_{13}^{2u_3} & \bmod \la c_{12},c_{13}^{2^\mu}\,\ra \\
    2 & c_{13}c_{12}^{v_1} & c_{12} & c_{12} & \bmod \la c_{13}^2\,\ra \\
    3 & c_{13} & 1 & 1 & \bmod \la c_{12},c_{13}^{2}\,\ra \\
    4 & c_{12}^{1+v_1} & 1 & c_{12} & \bmod \la c_{13}\,\ra \\
    5 & c_{12}^{v_1} & 1 & c_{12} & \bmod \la c_{13}\,\ra \\
    6 & c_{12}^{v_1} & 1 & c_{12} & \bmod \la c_{12}c_{13},c_{13}^{2}\,\ra \\
    7 & c_{13}c_{12}^{v_1} & 1 & c_{12} & \bmod \la c_{12}c_{13},c_{13}^{2}\,\ra
  \end{array} $$
\caption{}
\end{table}

\bigskip
For example, consider $j=1$. Then we may take $z=a_3$. Hence
$$t_1(\overline{a_1})\equiv a_1^2c_{13}\equiv c_{13}^{2u_1}c_{12}^{v_1}\equiv c_{13}^{2u_1}\bmod \la c_{12},c_{13}^{2^\mu}\ra.$$
A similar calculation works for $t_1(\overline{a_2}).$ Now for $a_3$ we have
$t_1(\overline{a_3})\equiv a_{3}^2\bmod H_1'$ and thus
$$t_1(\overline{a_3}^{2^{n-1}})\equiv a_3^{2^n}\equiv c_{13}^{2u_3}c_{12}\equiv c_{13}^{2u_3}\bmod \la c_{12},c_{13}^{2^\mu}\ra.$$
The rest of the entries are done in a similar fashion and are  left to the reader.

Once again, by comparing the class numbers of the fields $k_j$ with the indices  $(H_j:H_j')$ we see that $\{i_1,i_2\}=\{1,2\}$, i.e.,
$\{\Gal(k^2/k_1),\Gal(k^2/k_2)\}=\{H_1,H_2\}.$ Thus  exactly one of $\ker t_1$ or $\ker t_2$ has order $2$ while the other, as well as $\ker t_j$ for $j=3,\dots, 7$, has order $4$.  In particular, $\#\ker t_5=4$. But we can see that $\overline{a}_2\in \ker t_5$, whereas $\overline{a}_3^{2^{n-1}}\not\in \ker t_5$. This implies that
$$\ker t_5=\left\{
                \begin{array}{ll}
                 \la \overline{a}_1,\overline{a}_2\ra, & \hbox{if}\;\; v_1=0 \\
                  \la \overline{a}_1\overline{a}_3^{2^{n-1}},\overline{a}_2\ra, & \hbox{if}\;\; v_1=1.
                \end{array}
              \right.$$
If $v_1=1$, then we may change  variables $a_1'=a_1a_3^{2^{n-1}}, a_2'=a_2, a_3'=a_3$. This yields, arguing as above using (\ref{E5}), $$c_{12}'=[a_1',a_2']=c_{12}c_{13}^{2^{n}\tau}\,,\;\;c_{13}'=c_{13}c_{13}^{2^n\sigma}\,,\;\;c_{23}'=c_{23},$$
for some integers $\tau,\sigma$; and by (\ref{E6}) and   (\ref{E5})
$$a_1'^2c_{13}'=a_1^2a_3^{2^n}[a_3^{2^{n-1}},a_1][a_3^{2^{n-1}},a_1,a_3^{2^{n-1}}]c_{13}'=c_{13}^{2u_1}c_{12}c_{13}^{2u_3}c_{12}c_{13}^{2^nx}, $$
for some integer $x$. Therefore, $a_1'^2c_{13}'=c_{13}'^{2u_1'}$ for some integer $u_1'$. Also, we have
$$a_2'^2=c_{13}'^{2u_2'}\,,\;\;a_3'^{2^n}=c_{13}'^{2u_3'}c_{12}',$$
for some integers $u_2',u_3'$ (since $c_{23}\in G_3=\la c_{13}'^2\ra$ by (\ref{E4})\;). We claim that all the conditions in (\ref{E4}) are satisfied with respect to the new generators $a_j'$. The only one that needs to still be  verified is that $c_{12}'^2\in \la c_{13}'^{2^\mu}\ra.$  To see this, we have $H_1=\la a_1',a_2',a_3'^2, c_{12}',c_{13}'\ra,$ and
$H_1'\supseteq \la c_{12}',c_{13}'^{2^\mu}\ra.$ Since $(G_2:H_1')=2^\mu$,  we see that $c_{12}'^2$ must be in $\la c_{13}'^{2^\mu}\ra$.
Therefore, we may assume $v_1=0$.
\bigskip

We thus now have, without loss of generality,

\begin{equation}\label{E7}
a_1^2c_{13}=c_{13}^{2u_1}\,,\;\;a_2^2=c_{13}^{2u_2}\,,\;\; a_3^{2^n}=c_{13}^{2u_3}c_{12}.
\end{equation}

\bigskip

From the above we can conclude that
\begin{eqnarray*}
  \ker t_2 &=& \la \overline{a}_2\overline{a}_3^{2^{n-1}}\ra \\
  \ker t_3 &=& \la \overline{a}_2,\overline{a}_3^{2^{n-1}}\ra \\
  \ker t_4 &=& \la \overline{a}_2,\overline{a}_1\overline{a}_3^{2^{n-1}}\ra \\
 \ker t_5 &=& \la \overline{a}_1,\overline{a}_2\ra \\
 \ker t_6 &=& \la \overline{a}_1,\overline{a}_2\ra \\
 \ker t_7 &=& \la \overline{a}_2,\overline{a}_1\overline{a}_3^{2^{n-1}}\ra
\end{eqnarray*}
We do not have enough information yet to determine $\ker t_1$, but we do at least know that
$$t_1(\overline{a}_1)=c_{13}^{2u_1}H_1',\;\; t_1(\overline{a}_2)=c_{13}^{2u_2}c_{23}H_1',\;\;t_1(\overline{a}_3^{2^{n-1}})=c_{13}^{2u_3}H_1',$$
where $H_1'=\la c_{12},c_{13}^{2^\mu}\ra$.

Comparing data between the fields $k_j$ and the groups $H_j$, we see that $\Gal(k^2/k_j)=H_j$ for $j=1,2,3$.
Since $\Cl_2(k)\simeq \overline{G}$  via the Artin map, and $\overline{G}^2\cap \overline{G}[2]=\la \overline{a_3}^{2^{n-1}}\ra$, it follows  by (\ref{E1}) that  $[\mathtt{2}\fp]\leftrightarrow \overline{a}_3^{2^{n-1}}$, from above that $[\mathtt{\fp}]\leftrightarrow \overline{a}_2\overline{a}_3^{2^{n-1}}$, and thus that $[\mathtt{2}]\leftrightarrow\overline{a}_2.$ But then  since $[\mathtt{2}]\not\in\kappa_1$, $t_1(\overline{a}_2)\not=1$ (but clearly $t_1(\overline{a}_2^2)=1$),  by Table~4  $c_{13}^{2u_2}c_{23}\not\in \la c_{12},c_{13}^{2^\mu}\ra$, but $c_{13}^{4u_2}c_{23}^2\in \la c_{12},c_{13}^{2^\mu}\ra$;  therefore
$$c_{13}^{2u_2}c_{23}=c_{13}^{2^{\mu-1}u_0},\; \text{for some odd} \;u_0.$$
Similarly, since $[\mathtt{2}\fp]\not\in\kappa_1$, and so $t_1(\overline{a}_3^{2^{n-1}})\not=1$, it follows that
$$c_{13}^{2u_3}=c_{13}^{2^{\mu-1}v_0},\;\text{for some odd}\; v_0.$$
We also need some information on $u_1$. We don't know if $\overline{a}_1$ lies in $\ker t_1$; however, its square certainly does. Thus, for some integer $w$,
$$c_{13}^{2u_1}=c_{13}^{2^{\mu-1}w}.$$
Summarizing, we have
\begin{equation}\label{E8}
    c_{13}^{2u_2}c_{23}=c_{13}^{2^{\mu-1}u_0},\;\;c_{13}^{2u_3}=c_{13}^{2^{\mu-1}v_0},\;\;c_{13}^{2u_1}=c_{13}^{2^{\mu-1}w},
\end{equation}
with integers $u_0,v_0,w$ such that $u_0,v_0$ are odd.

\bigskip

To obtain more information on the group $G$, we will use the following results.
\begin{prop}\label{P2} Let $H_\gen=G^2=\la a_3^2,c_{13},c_{12}\ra.$ Then
$$H_\gen'=\la c_{13}^{2^\mu}\,\ra.$$
\end{prop}
\begin{proof} Since $H_\gen=\Gal(k^2/k_\gen)$, it follows that  by Proposition~\ref{P1},
$$2^{n+\mu}=h_2(k_\gen)=(H_\gen:H_\gen').$$
Since $(H_\gen:G_2)=2^{n-1}$, $(G_2:G_3)=4$, and $H_\gen'\subseteq G_3$, we have $(G_3:H_\gen')=2^{\mu-1}.$ Therefore since $G_3=\la c_{13}^2\ra$, we must have $H_\gen'=\la c_{13}^{2^\mu}\,\ra.$
\end{proof}

This leads to the following important proposition:
\begin{prop}\label{P2.5} Let $H_\gen$ be as given in Proposition~\ref{P2} and let $2^m$ be the order of $c_{13}$. Then $\mu=m$, and thus $H_\gen'=1$.
\end{prop}
\begin{proof}  We claim that $H_\gen'\subseteq \la c_{13}^{2^{\mu+1}}\,\ra.$ Now if $\mu\not=m$, in which case $\mu<m$  by (\ref{E1.5}), then our claim yields a  contradiction to Proposition~\ref{P2}. To establish the claim,  first notice that $H_\gen'=\la [c_{13},a_3^2],[c_{12},a_3^2], (H_\gen)_3\ra,$ and so we start by  showing that
$$[c_{13},a_3^2]\equiv [c_{12},a_3^2]\equiv 1\bmod \la c_{13}^{2^{\mu+1}}\,\ra.$$

We have $[c_{13},a_3^2]=c_{133}^2c_{1333}$. Now, (\ref{E4}) implies that
$$c_{133}=c_{13}^{2x_0},\;\;\text{for some odd}\;x_0.$$
By (\ref{E7}),
$1=[a_1^2,a_1]=[c_{13}^{-1+2u_1},a_1]=c_{131}^{-1+2u_1}$, which yields $c_{131}=1$ (since the exponent is odd). Moreover
by (\ref{E8}), $c_{13}^2=c_{13}^2c_{131}=[a_1^2,a_3]=[c_{13}^{-1+2^{\mu-1}w},a_3]=c_{133}^{-1+2^{\mu-1}w}.$ Therefore,
$$c_{13}^2c_{133}=c_{133}^{2^{\mu-1}w}=c_{13}^{2^\mu x_0w}.$$
But then
$$[c_{13},a_3^2]=c_{133}^2c_{1333}=[c_{13}^2c_{133},a_3]=c_{133}^{2^\mu x_0w}=c_{13}^{2^{\mu+1}x_0^2w}\in \la c_{13}^{2^{\mu+1}}\,\ra.$$

On the other hand, $[c_{12},a_3^2]=c_{123}^2c_{1233}.$ By (\ref{E7}), (\ref{E8}), and the previous paragraph,
$$1=[a_3^{2^n},a_3]=c_{133}^{2^{\mu-1}v_0}c_{123}\;\Rightarrow\;c_{123}=c_{133}^{-2^{\mu-1}v_0}=c_{13}^{-2^{\mu}v_0x_0},$$
and so $c_{123}^2,c_{1233}\in\la c_{13}^{2^{\mu+1}}\,\ra.$ Therefore $[c_{12},a_3^2]\in\la c_{13}^{2^{\mu+1}}\,\ra$. But then it follows that  $(H_\gen)_3\subseteq \la c_{13}^{2^{\mu+2}}\ra.$ Hence $H_\gen'\subseteq \la c_{13}^{2^{\mu+1}}\ra,$ as desired.
\end{proof}

With all this, we get
\begin{eqnarray}\label{E9}
c_{13}^{2u_1}=c_{13}^{2^{m-1}w}\,,\;\;c_{13}^{2u_2}=c_{23}c_{13}^{2^{m-1}}\,,\;\;c_{13}^{2u_3}=c_{13}^{2^{m-1}}\,,\\
\nonumber \Rightarrow\;a_1^2c_{13}=c_{13}^{2^{m-1}w}\,,\;\;a_2^2=c_{23}c_{13}^{2^{m-1}}\,,\;\;a_3^{2^n}=c_{12}c_{13}^{2^{m-1}}\,.
\end{eqnarray}

Also notice that by this proposition and (\ref{E4}), $c_{12}^2=1.$
Hence $H_1'=\la c_{12}\ra.$

\bigskip

We can simplify the relations yet a little more. Namely, we  will see that by a change of generators we can assume that $c_{23}=1$. To show this, first notice that since $c_{23}\in G_3=\la c_{13}^2\ra$,  we have $c_{23}=c_{13}^{2y}$, for some integer $y$, and thus $a_2^2=c_{13}^{2y}c_{13}^{2^{m-1}}.$ Now let
$$a_1'=a_1,\quad a_2'=a_2c_{13}^y,\quad a_3'=a_3.$$
Hence $c_{12}'=c_{12}[a_1,c_{13}^y]=c_{12}c_{131}^{-y}$ (by induction on $y$), and since $c_{131}=1$, as we saw in the proof of Proposition~\ref{P2.5}, we get $c_{12}'=c_{12}.$  Clearly $c_{13}'=c_{13}$. Moreover, $c_{23}'=c_{23}[c_{13}^y,a_3]=c_{23}c_{133}^y$ (again by induction on $y$), and since $c_{133}=c_{13}^{-2}$, by the proof of Proposition~\ref{P2.5}, we see that $c_{23}'=c_{23}c_{13}^{-2y}=1.$

Next, we have
$$a_1'^2=a_1^2=c_{13}'^{-1+2^{m-1}w}, \quad a_3'^{2^n}=a_3^{2^n}=c_{12}'c_{13}'^{2^{m-1}},$$
$$a_2'^2=a_2^2c_{13}^{2y}[c_{13}^y,a_2]=c_{13}^{2y}c_{13}^{2^{m-1}}c_{13}^{2y}[c_{13}^y,a_2]=c_{13}^{4y+2^{m-1}}[c_{13}^y,a_2].$$
But we claim that $[c_{13}^y,a_2]=1.$ To  show this, it suffices to prove that $c_{132}=1$, for then the claim follows by induction on $y$ again.   To this end, first we have
$$c_{12}^2c_{121}=[a_1^2,a_2]=[c_{13}^{-1+2^{m-1}w},a_2]=c_{132}^{-1}.$$
Since $c_{12}^2=1$, we get $c_{132}=c_{121}^{-1}.$  Moreover, from above we have $[a_1,a_3^2]=c_{13}^2c_{133}=1$, and so (by induction) $[a_1,a_3^{2^n}]=1$. Whence
$$1=[a_3^{2^n},a_1]=[c_{12}c_{13}^{2^{m-1}},a_1]=c_{121}.$$
Therefore $c_{132}=1$. This implies that $a_2'^2=c_{13}'^{4y+2^{m-1}}.$

From all this, we may replace (\ref{E9}) by
\begin{eqnarray*}
a_1^2=c_{13}^{-1+2^{m-1}w}\;(w\in\{0,1\}),\;\;a_2^2=c_{13}^{4y+2^{m-1}},\;\;a_3^{2^n}=c_{12}c_{13}^{2^{m-1}},\\
\nonumber c_{23}=c_{12}^2=c_{13}^{2^m}=1.\hspace{2in} \phantom{.}
\end{eqnarray*}

However, we can simplify $a_2^2$ quite a bit. Namely, write $a_2^2=c_{13}^x$. Then since $c_{23}=1$, we have $$1=[a_2^2,a_3]=[c_{13}^x,a_3]=c_{133}^x=c_{13}^{-2x}.$$
Hence $a_2^2=c_{13}^{2^{m-1}\varepsilon}$ for some $\varepsilon\in\{0,1\}.$

Finally, we claim that $w=0$ and $w=1$ produce isomorphic groups. For consider the change of generators given by
$$a_1'=a_1a_3^{2^{n-1}},\; a_2'=a_2,\;a_3'=a_3a_2.$$
Then it can be seen that
$$c_{12}'=c_{12},\; c_{23}'=c_{23}=1,\;c_{13}'=c_{12}c_{13},$$
and therefore
$$a_1'^2=a_1^2a_3^{2^n}=c_{13}^{-1+2^{m-1}w}c_{12}c_{13}^{2^{m-1}}=c_{13}'^{-1+2^{m-1}(1+w)}, $$
along with
$$a_2'^2=c_{13}'^{2^{m-1}\varepsilon},\;a_3'^{2^n}=c_{12}'c_{13}'^{2^{m-1}},\;c_{23}'=c_{12}'^2=c_{13}'^{2^m}=1.$$
Hence without loss of generality, we may take $w=0$.

Therefore we may assume that

$$a_1^2=c_{13}^{-1},\;a_2^2=c_{13}^{2^{m-1}\varepsilon},\;a_3^{2^n}=c_{12}c_{13}^{2^{m-1}},\;c_{23}=c_{12}^2=c_{13}^{2^m}=1.$$

\bigskip

Let
$$\Gamma_{n,m,\varepsilon}=\la a_1,a_2,a_3: a_1^2=c_{13}^{-1},\;a_2^2=c_{13}^{2^{m-1}\varepsilon},\;a_3^{2^n}=c_{12}c_{13}^{2^{m-1}},\;c_{23}=c_{12}^2=c_{13}^{2^m}=1\ra,$$
for $n,m>1$ and $\varepsilon=0,1$.
We now will show that if $\Gamma_{n,m,\varepsilon}=\Gal(k^\infty/k)$ for a field $k$ satisfying the conditions of Theorem~\ref{T1}, then $\varepsilon=1$. This will be established by computing the order of the kernel of the transfer map $t_{H_2,H_1\cap H_2}$ (cf.\,Lemma~\ref{L5} below) and comparing this via the class field correspondence  to the order of the capitulation kernel in the extension $k(\sqrt{p},i)/k$, as determined in Proposition~\ref{P2.75} below.
\begin{lem}\label{L5} Let $G=\Gamma_{n,m,\varepsilon}$ be as described above. Let $H_1$ and $H_2$ be the maximal subgroups of $G$ as given in Table~3. Let $t=t_{H_2,A}$ be the transfer map from $H_2/H_2'$ to $A/A'$ where $A=H_1\cap H_2$. Then
$$\#\ker t=\left\{
            \begin{array}{ll}
              8, & \hbox{if}\;\; \varepsilon=0, \\
              4, & \hbox{if}\;\;\varepsilon=1.
            \end{array}
          \right. $$
\end{lem}
\begin{proof} Recall from Table~3 that $H_1=\la a_1,a_2,a_3^2,c_{12},c_{13}\ra$ and $H_2=\la a_2,a_3,c_{12},c_{13}\ra$; thus $A=H_1\cap H_2=\la a_2,a_3^2,c_{12},c_{13}\ra$. By Table~4 we have $H_2'=\la c_{13}^2\ra.$ Now let $\overline{H}_2=H_2/H_2'$ and for $x\in H_2$ let $\overline{x}=xH_2'$. Then notice that
$$\overline{H}_2=\la \overline{a_2},\overline{a_3},\overline{c_{13}}\ra\simeq (2,2,2^{n+1}),$$
where $ \overline{a_2},\overline{a_3},\overline{c_{13}}$ have orders $2,2^{n+1},2$, respectively. Hence \\$\ker t\subseteq \overline{H}_2[2]=\la \overline{a_2},\overline{a_3}^{2^n},\overline{c_{13}}\ra$, which is of order $8$.

Now we compute $\ker t$. First observe that $A'=1$ (here we use previous relations as well as the fact that $c_{122}=1$ as the reader can verify) and that $H_2=A\cup Aa_3$.  Hence for $x\in H_2$,
$$t(\overline{x})=\left\{
                    \begin{array}{ll}
                      x^2[x,a_3], & \hbox{if}\;\; x\in A, \\
                      x^2, & \hbox{if not.}
                    \end{array}
                  \right.$$
We thus have
$t(\overline{a_2})=a_2^2c_{23}=c_{13}^{2^{m-1}\varepsilon},$ \;$t(\overline{a_3})=a_3^2,$ implying that $t(\overline{a_3}^{2^n})=a_3^{2^{n+1}}=1,$ and $t(\overline{c_{13}})=c_{13}^2c_{133}=1$. Therefore
$$\ker t=\left\{
           \begin{array}{ll}
             \la \overline{a_2},\overline{a_3}^{2^n},\overline{c_{13}}\ra , & \hbox{if}\;\;\varepsilon=0, \\
             \la \overline{a_3}^{2^n},\overline{c_{13}}\ra, & \hbox{if}\;\;\varepsilon=1.
           \end{array}
         \right.$$
This establishes the lemma.
\end{proof}

\begin{prop}\label{P2.75} Let $k$ be an imaginary quadratic number field as given in\\ Theorem~\ref{T1}. Let $k_2=k(\sqrt{p}\,)$ (as in Table~1) and let $K=k(\sqrt{p},i)=\Q(i,\sqrt{p},\sqrt{qq'}\,).$ Then the capitulation kernel $\kappa_{K/k_2}$ of the extension $K/k_2$ has order $4$.
\end{prop}
\begin{proof} First recall that the order of $\kappa=\kappa_{K/k_2}$ is given by $$\#\kappa=[K:k_2](E_{k_2}:N_{K/k_2} E_K)=2(\la -1,\varepsilon_p\ra:N_{K/k_2}E_K)$$ (see Table~1). We therefore need to determine $E_K$. We use Kuroda's class number formula on the $V_4$-extension $K/k$ (see \cite{Lem1}):
$$h_2(K)=\frac1{4}q(K/k)h_2(k_1)h_2(k_2)h_2(k_3)/h_2(k)^2,$$
where $k_j$ are as given in Table~1 and $q(K/k)=(E_K:E_1E_2E_3)$ with $E_j=E_{k_j}$. Hence, again by Table~1,
$$E_1E_2E_3=\la i,\varepsilon_p,\varepsilon_{qq'},\varepsilon_{pqq'}\ra.$$
Now $K=(k^2)^A$, the fixed field in $k^2$ of the subgroup $A$ given in Lemma~\ref{L5}. By class field theory $h_2(K)=(A:A')=\#A=2^{m+n+1}.$ Finally, given the values of $h_2(k_j)$ in Table~1, we see that $q(K/k)=2$.

We claim that $\sqrt{\varepsilon_{pqq'}}\in K$ which thus implies that $$E_K=\la i,\varepsilon_p,\varepsilon_{qq'},\sqrt{\varepsilon_{pqq'}}\ra.$$

To establish the claim, we will show that for $\eta=\sqrt{\varepsilon_{pqq'}}$,
$$\Q(\eta)=L, \;\text{with}\;\;L=\Q(\sqrt{p},\sqrt{qq'}\,).$$
(Since $L\subseteq K$, the claim will thus follow.)  We use Kuroda's class number formula on the $V_4$-extension $L/Q$:
$$h_2(L)=\frac1{4}q(L/\Q)h_2(p)h_2(qq')h_2(pqq'),$$
where $q(L/\Q)=(E_L:\la -1,\varepsilon_p,\varepsilon_{qq'},\varepsilon_{pqq'}\ra).$ Now $h_2(p)=h_2(qq')=1,$ and $h_2(pqq')=2.$  Since $L/\Q(\sqrt{pqq'}\,)$ is unramified and $h_2(pqq')=2,$ we see that $h_2(L)=1$. Thus $q(L/\Q)=2$. By considering the genus characters associated with the prime discriminants of $pqq'$, namely $\chi_p,\chi_{-q},\chi_{-q'}$, we see that $\chi(\fq\fq')=+1$ for each of these characters, where $\fq,\fq'$ are the prime ideals of $\Q(\sqrt{pqq'}\,)$ containing $q,q'$, respectively. Hence $\Q(\eta)=\Q(\sqrt{pqq'},\sqrt{qq'}\,)=L$, as desired.

Now we need to compute the norm $N=N_{K/k_2}$ of $E_K=\la i,\varepsilon_p,\varepsilon_{qq'},\eta\ra.$ To this end, we see that $Ni=1, N\varepsilon_p=\varepsilon_p^2,$ and $N\varepsilon_{qq'}=1.$ Hence we are left to compute $N\eta.$ Since $N\eta^2=N\varepsilon_{pqq'}=1$, we see that  $N\eta=\pm 1.$

We claim that $N\eta=-1$. To establish this, first we show  that $\eta=x\sqrt{p}+y\sqrt{qq'}$ for some $x,y\in \Q$. This can be seen as follows: since $\eta\in L=\Q(\sqrt{p},\sqrt{qq'}\,)$,  $\eta=w+x\sqrt{p}+y\sqrt{qq'}+z\sqrt{pqq'}$, for $w,x,y,z\in\Q$. Let $\tau\in\Gal(L/\Q),$  determined by  $\tau(\sqrt{p})=-\sqrt{p},\; \tau(\sqrt{pqq'})=\sqrt{pqq'}.$ Then we have
$$-\eta^2=\eta\tau(\eta),\;\;\text{and thus}\; -\eta=\tau(\eta),$$
which in turn implies that $w=z=0$. Hence since $\eta=x\sqrt{p}+y\sqrt{qq'}$, we see that $N\eta=N_{L/{\Q(\sqrt{p})}}\eta=px^2-qq'y^2=\pm 1.$ Since $(p/q)=-1$, we have $N\eta=-1$, as desired.

Thus $N_{K/k_2}(E_K)=\la -1,\varepsilon_p^2\ra$ and so $(E_{k_2}:N_{K/k_2}E_K)=2$. Therefore $\#\kappa=4$, as desired.
\end{proof}

By the class field correspondence given by the Artin map, we have $\#\ker t=\#\kappa=4$, where $t$ is as given in Lemma~\ref{L5} and $\kappa$ as in Proposition~\ref{P2.75}. Therefore, $\varepsilon=1,$ i.e.,
$$\Gamma_{n,m,\varepsilon}=\Gamma_{n,m,1}=\Gamma_{n,m}.$$

Conversely, if $k$ is an imaginary quadratic number field for which \\$G=\Gal(k^\infty/k)=\Gamma_{n,m}$, for some $n,m>1$, then it follows that $G/G'\simeq (2,2,2^n)$, $G'\simeq (2,2^m)$, and $G/G^{(2,2)}=32.012$. By \cite{BS}, the fields of Type 4(p) (as given in \cite{BS}) are the only ones for which this can happen. These are precisely the fields whose discriminants satisfy (\ref{ET}).

All of this completes the proof of the theorem.

\bigskip
\section{\bf Some Consequences of the Theorem}

As a consequence of Theorem~\ref{T1} we obtain the following corollary concerning the $2$-class groups of some of the intermediate fields in the extension $k^1/k$. Among the fields we consider, perhaps the most interesting is the unramified quadratic field $\Q(\sqrt{-p},\sqrt{qq'}\,)$, as the structure of its $2$-class group depends on the relative sizes of $n$ and $m$; cf.\,the next corollary.
\begin{cor}\label{C2} Let $k$ be a  quadratic number field as given in Theorem\,\ref{T1}, $2^n=h_2(k)/4$, and $2^m=h_2(-p)$.
Then the following statements are valid:

If as in Table\,1, $k_1=\Q(\sqrt{-p},\sqrt{qq'}\,)$, then
$$\Cl_2(k_1)\simeq \left\{
                     \begin{array}{ll}
                       (2,2^{m+1},2^{n-1}), & \hbox{if}\;\;m\geq n-1, \\
                       (2,2^{m},2^{n}), & \hbox{if}\;\;m < n-1.
                     \end{array}
                   \right.$$
Moreover,
\begin{eqnarray*}
  \Cl_2(k(\sqrt{p}\,)) &\simeq & (2,2,2^{n+1}),\\
  \Cl_2(k(i)) &\simeq & (2,2,2^{n}),\\
  \Cl_2(K) &\simeq & (2,2^{n+1}),\\
  \Cl_2(k_\gen) & \simeq & (2^n,2^m),\\
  \Cl_2(k^1) & \simeq & (2,2^m),
\end{eqnarray*}
where $K=k_j$ for $j=4,5,6,7$ in Table 1, i.e.,
$$K\in\{k(\sqrt{-q}\,),k(\sqrt{-q'}\,),k(\sqrt{q}\,),k(\sqrt{q'}\,)\}.$$
\end{cor}
\begin{proof}(Sketch) We consider the results in reverse order. Notice that $\Cl_2(k^1)\simeq (2,2^m)$ follows directly from Theorem~\ref{T1}.

Next, consider $k_\gen=\Q(i,\sqrt{p},\sqrt{q},\sqrt{q'}\,)$. As we saw in the proof of \\Proposition~\ref{P2},
$\Gal(k^2/k_\gen)=H_\gen=\la a_3^2,c_{12},c_{13}\ra$, and by Proposition~\ref{P2.5}, $H_\gen'=1$. Thus $\Cl_2(k_\gen)=H_\gen$.
Now, by the presentation of $G$ in the theorem, we have $c_{12}=a_3^{2^n}c_{13}^{2^{m-1}}$ and so $H_\gen=\la a_3^2,c_{13}\ra.$ Moreover, $\la a_3^2\ra\cap\la c_{13}\ra=1$, $a_3^2$ has order $2^n$, and $c_{13}$ has order $2^m$. Thus
$\Cl_2(k_\gen)=H_\gen\simeq (2^n,2^m).$

Concerning the rest of the results, by comparing the capitulation kernels $\kappa_j$ in Table~1 with the kernels $\ker t_j$, we can see that $H_j=\Gal(k^2/k_j)$ for $j=1,2,3$. So let $H_4=\la a_1,a_3,c_{12},c_{13}\,\ra,$ in which case $H_4'=\la c_{13}\ra$ (cf.\,Table~3 above). Let $\overline{H_4}=H_4/H_4'$ and $\overline{h}=hH_4'$ for any $h\in H_4$. Hence $\overline{H_4}\simeq \Cl_2(K)$, for some $K\in \{k_j:j=4,\dots, 7\}.$ By the presentation of $G$ in the theorem, we have
$\overline{H_4}=\la \overline{a_1},\overline{a_3}\,\ra$, since $\overline{c_{12}}=\overline{a_3}^{2^n}$ and $\overline{c_{13}}=1$. But  $\overline{a_1}^2=1$ and $\overline{a_3}^{2^{n+1}}=1$. Therefore, $\Cl_2(K)\simeq \overline{H_4}\simeq (2,2^{n+1}).$  For $j=5,6,7$,\; $H_j$ are handled in a similar manner and left to the reader.

Now consider $k_3=k(i)$ as in Table~1. Then $\Gal(k^2/k_3)=H_3=\la a_1a_3,a_2,c_{12},c_{13}\ra$ and so by Table~3, $H_3'=\la c_{12},c_{13}^2\ra.$ Let $\overline{H_3}=H_3/H_3'\simeq \Cl_2(k_3)$. Then $\overline{H_3}=\la \overline{a_1}\,\overline{a_3}, \overline{a_2},\overline{c_{13}}\ra$, since $\overline{c_{12}}=1$. But $(\overline{a_1}\,\overline{a_3})^{2^n}=1$, $\overline{a_2}^2=\overline{c_{13}}^2=1$. Thus $\Cl_2(k_3)\simeq \overline{H_3}\simeq (2,2,2^{n}).$   The case of $k_2$ is similar and again left to the reader.

Finally consider $k_1=k(\sqrt{-p}\,)$. Hence $\Gal(k^2/k_1)=H_1=\la a_1,a_2,a_3^2,c_{12},c_{13}\ra$ and $H_1'=\la c_{12}\ra.$ Let $\overline{H_1}=H_1/H_1'$; then $\overline{H_1}=\la \overline{a_1},\overline{a_2},\overline{a_3}^2\ra$, since $\overline{c_{12}}=1$ and $\overline{c_{13}}=\overline{a_1}^2.$  We now have
$$\overline{a_1}^{2^m}=\overline{a_2}^2=(\overline{a_3}^2)^{2^{n-1}}=\overline{a_3}^{2^n}=\overline{c_{13}}^{2^{m-1}}.$$
Let $\beta=\overline{a_2}\,\overline{a_1}^{2^{m-1}}$; then $\overline{H_1}=\la \beta,\overline{a_1},\overline{a_3}^2\ra$ and $\beta^2=1$. Notice then that $\la \beta\ra\cap \la \overline{a_1},\overline{a_3}^2\ra=1$ (cf.\,Burnside Basis Theorem), whence $\overline{H_1}\simeq \la \beta\ra \oplus \la \overline{a_1},\overline{a_3}^2\ra.$

 From above we have $\overline{a_1}^{2^m}=\overline{a_3}^{2^n}$, each of order $2$. Hence $\la \overline{a_1}\ra\cap \la \overline{a_3}^2\ra=\la\overline{a_1}^{2^m}\ra=\la \overline{a_3}^{2^n}\ra$.   We now consider two cases: $m<n-1$ and $m\geq n-1$.

Case 1. $m<n-1$.  Let $\alpha=\overline{a_1}\,\overline{a_3}^{\,-2^{n-m}}$. Then $\alpha^{2^m}=1$, $\overline{a_3}^{2^{n+1}}=1$ (where the exponents are minimal), and $\la \alpha\ra\cap \la \overline{a_3}^2\ra=1.$ Thus $\la \overline{a_1},\overline{a_3}^2\ra=\la \alpha,\overline{a_3}^2\ra \simeq (2^m,2^n)$ and so $\overline{H_1}\simeq (2,2^m,2^n).$

Case 2. $m\geq n-1$. This time let $\alpha=\overline{a_3}^2\,\overline{a_1}^{\,-2^{m+1-n}}$; then
the order of $\overline{a_1}$ is $2^{m+1}$ and that of $\alpha$ is $2^{n-1}$. Arguing as above, it follows that
$\overline{H_1}\simeq (2,2^{m+1},2^{n-1}).$
\end{proof}

\bigskip

\section{\bf  Some Extensions of Theorem~\ref{T1}\label{S7}}

For the fields $k$ given in Theorem~\ref{T1} we found candidates for
$\Gal(k^2/k)$ as
$$\Gamma_{n,m,\varepsilon}=\la a_1,a_2,a_3: a_1^2=c_{13}^{-1},\;
a_2^2=c_{13}^{2^{m-1}\varepsilon},\;a_3^{2^n}=c_{12}c_{13}^{2^{m-1}},\;
c_{23}=c_{12}^2=c_{13}^{2^m}=1\ra,$$
for $n,m>1$ and $\varepsilon=0,1$. The question then naturally arises as
to what happens when $n=1$ or $m=1$. We consider those cases now, i.e.,

(I) $\Gamma_{1,1,\varepsilon}$, \;\;
(II) $\Gamma_{1,m,\varepsilon},\;m>1,$ \;\;
(III) $\Gamma_{n,1,\varepsilon},\;n>1.$

\begin{prop}\label{PS1} Let $k$ be an imaginary quadratic number field.

\medskip
\noindent (i) $\Gamma_{1,1,1}\simeq \Gal(k^\infty/k)$
\medskip

if and only if

\noindent the discriminant of $k$ factors into prime discriminants $d_j$, divisible by the prime $p_j$, such that  $d_k=d_1d_2d_3d_4$,

\noindent
1. $d_k\not\equiv 4\bmod 8$,

\noindent
2. $d_j<0$ for $j=1,2,3$ and $d_4>0$  (without loss of generality),

\noindent
3. $(d_j/p_4)=-1$ for $j=1,2,3$ and $(d_1/p_2)=(d_2/p_3)=(d_3/p_1)=-1$.

\bigskip

\noindent (ii) $\Gamma_{1,1,0}$ is not isomorphic to $\Gal(k^\infty/k)$ for any imaginary quadratic number field $k$.
\end{prop}
\begin{proof} By using GAP (or even by hand) we found  that
$$\Gamma_{1,1,0}=[32,29]=32.037,\quad \Gamma_{1,1,1}=[32,33]=32.041,$$
where the first term is the group ID as numbered  in GAP, the first coordinate being the order of the group; the latter term is as numbered in HS, i.e., Hall and Senior, \cite{HS}.  The first part of the proposition (i) then follows by Proposition~4 in \cite{BLS1}.
However, by Proposition~10 in \cite{BLS1}, there are no such fields $k$ with $\Gal(k^\infty/k)$ of type HS~$32.037$, since the fields there do not have elementary $\Cl_2(k^1)$, whereas $\Gamma_{1,1,0}'\simeq (2,2)$.
\end{proof}

\begin{prop}\label{PS2} There are no imaginary quadratic fields $k$ for which
$$\Gal(k^\infty/k)\simeq \Gamma_{1,m,\varepsilon},$$
for any $m>1$ and $\varepsilon=0,1$.
\end{prop}
\begin{proof} For brevity write $\Gamma=\Gamma_{1,m,\varepsilon}.$ Then $\Gamma/\Gamma_3=[32,30]=32.038$, by GAP (or by hand), and by Proposition~8 in \cite{BLS1} there exist imaginary quadratic fields $k$  for which $G=\Gal(k^\infty/k)$ satisfies $G/G_3$ of type $32.038$. But again by GAP, $G/G_4\simeq [64,166]$; however, $\Gamma/\Gamma_4$ is isomorphic to $[64,165]$ for $\varepsilon=0$ and  to $[64,163]$ when $\varepsilon=1$. Thus these $\Gamma$ are not realized as Galois groups for $k$ as above.
\end{proof}

\begin{prop}\label{PS3} Let $k$ be an imaginary quadratic number field.

\medskip
\noindent (i) $\Gamma_{n,1,1}\simeq \Gal(k^\infty/k)$ for some $n>1$
\medskip

if and only if

\noindent the discriminant $d_k=-4pqq'$, where $p,q,q'$ are primes such that
\begin{equation*}
 p\equiv 5\bmod 8,\;q,q'\equiv 3\bmod 8,\;  (p/q)=(p/q')=(-q/q')=-1.
\end{equation*}

\bigskip

\noindent (ii) $\Gamma_{n,1,0}$, for any $n>1$,  is not isomorphic to $\Gal(k^\infty/k)$ for any imaginary quadratic number field $k$.
\end{prop}

The proof of this result is somewhat more involved than the two previous propositions and will be established in the two lemmas below. But first we observe by GAP (or by hand) that  if we let $\Gamma=\Gamma_{n,1,\varepsilon}$, then
$$\Gamma/\Gamma^{(2,2)}\simeq [32,25]=32.014, \;\;\text{resp.}\;\; [32,26]=32.015,$$
for $\varepsilon=0$, resp. $1$.

First consider $\Gamma=\Gamma_{n,1,1}$. Since $\Gamma/\Gamma^{(2,2)}\simeq 32.015$, we see by \cite{BS} that $k$ is of Type 4(r), i.e., satisfies the conditions in Proposition~\ref{PS3}(i). But now we have the following lemma.
\begin{lem}\label{L6} Let $k$ be an imaginary quadratic number field with discriminant satisfying the conditions in Proposition~\ref{PS3}(i). Then
$$\Gal(k^\infty/k)\simeq \Gamma_{n,1,1},$$
where $2^n=h_2(k)/4.$
\end{lem}
\begin{proof}
We first claim that $\Gal(k^\infty/k)$ for $k$ of Type 4(r) can be presented  as follows:

$\Gal(k^\infty/k)\simeq \Gamma_n^{(4(r))},$ where $\Gamma_n^{(4(r))}$ is given by
$$\Gamma_n^{(4(r))}=\la a_1,a_2,a_3: a_1^2=a_2^2=c_{12},a_3^{2^n}=c_{13},c_{23}=c_{ij}^2=1,\; i,j=1,2,3\ra.$$

Let $G=\Gal(k^\infty/k).$ Then for starters we have $G=\la
a_1,a_2,a_3\ra$, where (without loss of generality)
$$a_1^2\equiv a_2^2\equiv a_3^{2^n}\equiv 1 \bmod G'\;\;(n>1).$$
Again by \cite{BS} we see that $H=G^{(2,2)}=\la a_1^2a_2^{-2},a_3^4,a_1^2c_{12},c_{23},c_{13}\ra$, since $G/H$ is of type 32.015. By Proposition~3 of \cite{BS} we have $H\cap G_2=\la a_3^{2^n}\ra,$ since $G_3=1.$

Now, we claim that we may assume $c_{23}=1$, without loss of
generality. For we have $G'=\la c_{12},c_{13},c_{23}\ra$, and
moreover, $c_{12}\not\in H$. Thus $H\cap G'=\la
c_{13},c_{23}\ra$. Hence (since $G'\simeq (2,2)$ by Theorem~7 in
\cite{BS}) $G'=\la c_{12},\gamma\ra$, where
$\gamma\in\{c_{13},c_{23},c_{13}c_{23}\}.$ If $c_{23}=1$, then we are
done.  If $c_{13}=1$, simply permute the indices $2$ and $1$. On the
other hand, if $c_{13}c_{23}=1$, then apply the substitutions
$a_1'=a_1,\,a_2'=a_2a_1,\,a_3'=a_3.$ Thus the claim is established.

Thus we may  assume that $G'=\la c_{12},c_{13}\ra,$ and hence $a_3^{2^n}=c_{13}.$ So far we have
$$a_1^2=c_{12}c_{13}^u,\,a_2^2=c_{12}c_{13}^v,\,a_3^{2^n}=c_{13},\,c_{23}=1,\,c_{ij}^2=1 \;\;(u,v\in\{0,1\}).$$

Let
$$\Gamma_n^{(u,v)}=\la a_1,a_2,a_3:a_1^2=c_{12}c_{13}^u,\,a_2^2=c_{12}c_{13}^v,\,a_3^{2^n}=c_{13},\,c_{23}=1,\,c_{ij}^2=1\ra,$$
where $u,v=0,1$. However, notice that $\Gamma_n^{(u,v)}=\Gamma_n^{(0,0)}$ by the substitutions: $a_1'=a_1a_3^{2^{n-1}u},\,
a_2'=a_2a_3^{2^{n-1}v},\,a_3'=a_3.$ We thus have $\Gamma_n^{(0,0)}\simeq \Gamma_n^{(4(r))}\simeq \Gal(k^\infty/k).$

Finally, we claim that $\Gamma_{n,1,1}=\Gamma_n^{(4(r))}.$ For
$$\Gamma_{n,1,1}=\la a_1,a_2,a_3: a_1^2=c_{13},\;a_2^2=c_{13},\;a_3^{2^n}=c_{12}c_{13},\;c_{23}=c_{12}^2=c_{13}^{2}=1\ra,$$
and changing the generators to $a_1'=a_1a_3^{2^{n-1}},$ $a_2'=a_2a_3^{2^{n-1}},$ and $a_3'=a_3a_2$ gives the presentation in
$\Gamma_n^{(4(r))}$.  This establishes the lemma and thus Proposition~\ref{PS3}(i).
\end{proof}

\begin{lem}\label{L7} Let $\Gamma=\Gamma_{n,1,0}$ with $n>1$. Then there is a maximal subgroup $H$ of $\Gamma$ such that the kernel of the transfer map, $t_{\Gamma,H}$, has order $8$.
\end{lem}
\begin{proof} We have
  $$\Gamma=\Gamma_{n,1,0}=\la a_1,a_2,a_3: a_1^2=c_{13},\;a_2^2=1,\;
  a_3^{2^n}=c_{12}c_{13},\;c_{23}=c_{12}^2=c_{13}^{2}=1\ra.  $$
  Now let $H=\la a_2a_3,a_1,c_{12}, c_{13}\ra$; hence $H'=\la
  c_{12}c_{13}\ra$ (here we have used the fact that $\Gamma_3=1$ as
  can easily be seen). Then we have the following values of the
  transfer map $t=t_{\Gamma,H}$:
  $$ t(\overline{a_1})=a_1^2c_{13}H'=1,\;t(\overline{a_2})=a_2^2H'=1,\;
     t(\overline{a_3}^{2^{n-1}})=c_{12}c_{13}H'=1,$$
  where $\overline{a_i}=a_i\Gamma'$.  Therefore, $\#\ker(t)=8$.     
\end{proof}

This now shows the validity of Proposition~\ref{PS3}(ii), since there
are no unramified quadratic extensions of an imaginary quadratic field
with capitulation kernel of order $>4$\;
(cf.\,\cite{Iwa},\,\cite{Che}).

\bigskip

\section{\bf  Examples}

We have used Pari and GAP for verifying our results in some special cases.
Observe first that, for $m + n \le 5$, the groups $\Gamma_{n,m,\eps}$
have the following ID numbers in GAP:
  $$ \begin{array}{rr|cc}
    n & m  & \eps = 0 & \eps = 1 \\ \hline
    2 & 2  & [ 128,  305 ] & [ 128, 306 ] \\
    3 & 2  & [ 256, 5004 ] & [ 256, 5006 ] \\
    2 & 3  & [ 256, 5023 ] & [ 256, 5024 ]
  \end{array} $$

We now give examples of imaginary quadratic number fields $k$ for which
$\Gal(k^2/k) = \Gamma_{n,m,1}$ for each of these three cases.

\section*{$G = [128, 306]$}
Let $k = \Q(\sqrt{-561}\,)$; here $p = 17$, $q = 3$ and $q' = 11$.
The seven unramified quadratic  extensions $K = k(\sqrt{m}\,)$ of $k$
have the following $2$-class groups:
$$ \begin{array}{c|ccccccc}
  m        & -1 & -q & q & -q' & q' & p & -p \\ \hline
  \Cl_2(K) & (2,2,4) & (2,8) & (2,8) & (2,8) & (2,8) & (2,2,8) & (2,2,8)
\end{array} $$
The square roots of $\alpha=\alpha_1 = 5 + 4 \sqrt{-17}$ and $\alpha_2
= 4 + \sqrt{33}$ generate\footnote{The explicit construction of
  unramified $2$-extensions is discussed in detail in \cite{LemT}; for
  the case of quaternion extensions see also \cite{LemQ}.} cyclic
quartic unramified $C_4$-extensions of $k$; the square root of $\beta
= (4i + \sqrt{17}\,)(4+\sqrt{17}\,)$ generates an unramified
quaternion extension, and the square roots of $\gamma=\gamma_1 = 1+4i$
and $\gamma_2 = \frac{1+\sqrt{17}}2$ generate nonnormal quartic
unramified extensions of $k$.

The unramified extensions $K_1 = k(i, \sqrt{\alpha_1}\,)$,
$K_2 = k(\sqrt{\beta}\,)$ and $K_3 = k(\sqrt{17}, \sqrt{\gamma_2}\,)$
have degree $8$ over $k$; their Galois groups are
$\Gal(K_1/k) \simeq C_2 \times C_4$,
$\Gal(K_2/k) \simeq H_8$ and $\Gal(K_3/k) \simeq D_4$. The
compositum $L = K_1K_2K_3$ is an unramified extension of relative degree $16$
over $k$ with Galois group $H_8 \curlywedge C_4$, the group HS~$16.010$ or
$[16,4]$ in GAP and pari. Similar remarks apply to the other cases
discussed below.

The field $M = k_\gen(\sqrt{\alpha}, \sqrt{\gamma}\,)$ has Galois group
$\Gal(M/\Q) = [64, 203]$.

\begin{figure}[ht!]
  \begin{center}
    \begin{tikzcd}[every arrow/.append style={dash}]
    &   & L \arrow[ld] \arrow[d] \arrow[rd]  &    &       \\
        & K_1 \arrow[d] \arrow[rd] \arrow[ld] & K_2 \arrow[d]
        & K_3 \arrow[d] \arrow[rd] \arrow[ld] &    \\
        k(\sqrt{\alpha}) \arrow[rd] & k(\sqrt{-\alpha}) \arrow[d]
        & {k(i,\sqrt{p})} \arrow[ld] \arrow[d] \arrow[rd]
        & k(\sqrt{\gamma_1}) \arrow[ld]  & k(\sqrt{\gamma_2}) \arrow[ld] \\
    & k(\sqrt{-p}) \arrow[rd]     & k(i) \arrow[d]
    & k(\sqrt{p}) \arrow[ld]     &           \\
    &  & k  & &                              
      \end{tikzcd}      
    \caption{Hasse diagram of $L/\Q$ with Galois group $H_8 \curlywedge C_4$}
  \end{center}
\end{figure}

\bigskip
Since $\Cl(k_\gen) = (4,4)$ and $\Cl(k^1) = (2,4)$, the $2$-class
field tower terminates at $k^2 = k_\gen^1$ (cf.\,Proposition~7 in
\cite{BLS}), and $\Gal(k^2/k)$ is a group of order $128$ with Schur
multiplier of order $2$, cf.\,\cite{Iwa} (also see \cite{Roq}
pp.\,235--236, where $H^{-3}(G,\Z)$ is the Schur multiplier of $G$).

Using GAP, we have verified that the only groups $G$ of order $128$
with $G' \simeq (2,4)$, $G/G' \simeq (2,2,4)$, Schur multiplier of order $2$
and subgroups of index $2$ with abelianizations corresponding to the
class groups in the table above are $[128, 305]$ and $[128, 306]$.

The abelianizations of the normal subgroups $H$ of index $\le 8$ in $G
= [128, 305] = \Gamma_{2,2,0}$ and $G = [128, 306]= \Gamma_{2,2,1}$
coincide. However, the abelianizations of the nonnormal subgroups $H$
of index $4$ in $G$ allow us to separate these two groups:
$$ \begin{array}{c|cc||c|cc}
   G        &  H & H/H' &  G        &  H & H/H' \\ \hline
  [128, 305] & [32,36]_2 & (2,2,8) & [128, 306]  & [32,37]_2 & (2,2,4) \\
             & [32,12]_2 & (2,8)   & & [32,12]_2 & (2,8) \\
             & [32, 4]_2 & (4,4)   & & [32, 3]_2 & (4,8) \\
             & [32, 5]_2 & (2,8)   & & [32, 5]_2 & (2,8)
  \end{array} $$
(Here, the subscripts $2$ indicate  the number of the corresponding subgroups.)
With $\gamma = \frac{1+\sqrt{17}}2$, the extensions $N/k$ corresponding
to some of these subgroups are the following:
$$ \begin{array}{l|l||l|l}
  N & \Cl_2(N) &  N & \Cl_2(N) \\ \hline
  k(\sqrt{\gamma}\,)     & (2,8)   & k(\sqrt{ \gamma'}\,)   & (4,8) \\
  k(\sqrt{p \gamma}\,)   & (2,8)   & k(\sqrt{-\gamma'}\,)   & (4,8) \\
  k(\sqrt{q \gamma}\,)   & (2,2,4) & k(\sqrt{ q\gamma'}\,)  & (2,8) \\
  k(\sqrt{q'\gamma}\,)   & (2,2,4) & k(\sqrt{-q\gamma'}\,)  & (2,8)
\end{array} $$
Thus $G = [128, 306]$.

In all cases we have considered, the two unramified quaternion extensions
$k(\sqrt{\beta}\,)$ and $k(\sqrt{q\beta}\,)$
of $k$ have $2$-class groups of type $(2,2,4)$ and $(4,4)$, although we
do not know beforehand which extension has which class group.

In addition, we have computed all imaginary quadratic number fields $k$
with discriminant $> -50,000$ and the following properties:
$\Cl_2(k) = (2,2,4)$ and $\Cl_2(k_\gen) = (4,4)$. Let $M$ denote the
maximal elementary abelian unramified $2$-extension of $k_\gen$; then
in all cases, $\Gal(M/\Q) \simeq [64,226]$ or $\simeq [64, 203]$.
Group theoretical calculations then show that in all of these cases,
$\Gal(k^2/k) \ne [128, 305]$, and that $\Gal(M/\Q) \simeq [64,203]$
implies $\Gal(k^2/k) \simeq [128, 306]$. The only such fields in the
range considered are the following:

$$ \begin{array}{r|rrr|ccc}
  d     & p & q & r & \alpha & \beta & \gamma \\ \hline
  - 2244 &  17 &   3 &  11 &  5 +  4 \sqrt{-17}
         & (4i+\sqrt{17}\,)(4+\sqrt{17}\,) &  1+4i \\
  -20292 &  89 &  19 &   3 &  1 +  8 \sqrt{-89}
         & (52i+\sqrt{89}\,)(8+\sqrt{89}\,) &  5+8i \\
  -26724 &  17 &   3 & 131 & 11 +  4 \sqrt{-17}
         & (52i + 7\sqrt{17}\,)(4+\sqrt{17}\,) &  1+4i \\
  -30756 & 233 &   3 &  11 & 43 +  4\sqrt{-233}
         & (8i+\sqrt{233}\,)(8+\sqrt{233}\,) & 13+8i \\
  -33252 &  17 & 163 &   3 &  7 + 16 \sqrt{-17}
         & (8i + 5\sqrt{17}\,)(4+\sqrt{17}\,) &  1+4i \\
  -46308 &  17 &   3 & 227 & 71 +  8 \sqrt{-17}
         & (16i + 5\sqrt{17}\,)(4+\sqrt{17}\,) &  1+4i
\end{array} $$

\bigskip

\section*{$G = [256, 5006]$}

The basic strategy in verifying the Galois group is the same as for
$G = [128, 306]$, so we will be brief.

The smallest example is $d = -4pqq'$ with $p = 17$, $q = 3$ and $q' = 107$.
We have $\Cl_2(k) = (2,2,8)$ and $\Cl_2(k_\gen) = (4,8)$; the seven quadratic
unramified extensions $K = k(\sqrt{m}\,)$ of $k$ have the following
$2$-class groups:
$$ \begin{array}{c|ccccccc}
  m        & -1 & -q & q & -q' & q' & p & -p \\ \hline
  \Cl_2(K) & (2,2,8) & (2,16) & (2,16) & (2,16) & (2,16) & (2,2,16) & (2,4,8)
\end{array} $$
The square roots of $\alpha = -7 +   4 \sqrt{-17}$,
$\beta = (4i + 13 \sqrt{17}\,)(4 + \sqrt{17}\,)$ and $\gamma = 1+4i$
generate unramified extensions of $k$ with Galois groups $C_4$ and $H_8$
for $\alpha$ and $\beta$, whereas $k(\sqrt{17}, \sqrt{\gamma}\,)$
is an unramified $D_4$-extension of $k$.  The $2$-class groups of the
two unramified quaternion extensions of $k$ are $\simeq (2,2,8)$ and
$(4,4)$.

Using GAP we have verified that the only groups $G$ of order $256$ with
the required properties are $[256, 5004]$ and $[256, 5006]$.
The $2$-class groups of the nonnormal unramified quadratic  extensions of $k$
allow us to conclude that $\Gal(k^2/k) = [256, 5006]$.

The fields with $d > -100\,000$ such that
\begin{itemize}
\item $\Cl_2(k) \simeq (2,2,8)$, $\Cl_2(k_\gen) = (4,8)$;
\item $\Gal(M/\Q) \simeq [64,203]$
\end{itemize}
are the following:
$$ \begin{array}{r|rrr|ccc}
  d & p & q & q' & \alpha & \beta & \gamma \\ \hline
  - 21828 & 17 & 3 & 107 & 7 + 4\sqrt{-17}
         & (4i + 13\sqrt{17}\,)(4+\sqrt{17}\,) & 1+4i \\
  - 28356 & 17 & 3 & 139
         & 59 + 4\sqrt{-17} & (20i + \sqrt{17}\,)(4+\sqrt{17}\,)   & 1+4i \\
  -91428 & 401 & 3 & 19
         & 119+2\sqrt{-401} & (32i + \sqrt{401}\,)(20 + \sqrt{401}\,) & 1+20i \\
  -97988 &  17 & 11 & 131
         & 109 + 8\sqrt{-17} & (8i + 9\sqrt{17}\,)(4 + \sqrt{17}\,) & 1+4i
\end{array} $$

\section*{$G = [256, 5024]$}

Consider the quadratic number field $k$ with discriminant
$d = -4 \cdot 3 \cdot 11 \cdot 41$. Using {\tt pari} we find that
$\Cl(k) = (2, 2, 4)$, and that the class groups of the
unramified quadratic extensions $K = k(\sqrt{m}\,)$ are given by:
$$ \begin{array}{c|ccccccc}
     m     &   -1    &  -3 & 3 & -11 & 11 & 41 & -41 \\ \hline
  \Cl_2(k) & (2,2,4) & (2,8) & (2, 8) & (2,8) & (2,8) & (2,2,8) & (2,2,16)
  \end{array} $$
Showing that $\Gal(k^2/k) \simeq [256,5023]$ or $[256,5024]$  is
straightforward (as above). The class groups of the nonnormal  unramified octic
extensions of $k$ allow us to conclude that $\Gal(k^2/k) \simeq [256,5024]$.

Fields with $d > -100,000$:

$$ \begin{array}{r|rrr|ccc}
      d  &   p & q & q' & \alpha & \beta & \gamma \\ \hline
  - 5412 &  41 &  3 &  11 & 13 +  4 \sqrt{-41}
         & (16i + \sqrt{41}\,)(4 + \sqrt{41}\,) &  5 + 4i \\
  -34276 &  41 & 11 &  19 & 35 +  4 \sqrt{-41}
         & (72i + \sqrt{41}\,)(4 + \sqrt{41}\,) &  5 + 4i\\
  -58308 & 113 &  3 &  43 & 37 + 20 \sqrt{-113}
         & (4i + \sqrt{113}\,)(8 + \sqrt{113}\,) &  7 + 8i \\
  -70692 & 137 &  3 &  43 & 41 +  8 \sqrt{-137}
         & (32i + \sqrt{137}\,)(4 +  + \sqrt{137}\,) & 11 + 4i  \\
  -88068 &  41 &  3 & 179 & 85 + 16  \sqrt{-41}
         & (92i + 11 \sqrt{41}\,)(4 + \sqrt{41}\,) & 5 + 4i \\
  -90852 & 113 &  3 &  67 &  1 +  4 \sqrt{-113}
         & (116i + 5 \sqrt{113}\,)(8 + \sqrt{113}\,) &  7 + 8i
\end{array} $$

\subsection*{Fields with Galois Group $\Gamma_{n,m,1}$}

The following table gives examples of imaginary quadratic number fields
for which $\Gal(k^2/k) \simeq \Gamma_{n,m}(=\Gamma_{n,m,1})$ and $n + m \le 8$:

$$ \begin{array}{r|rrr|rr}
        d   &  p   &   q &   q' &   m &  n \\ \hline
   -  2244  &   17 &   3 &  11  &   2 &  2 \\
   - 21828  &   17 &   3 & 107  &   2 &  3 \\
   -  5412  &   41 &   3 &  11  &   3 &  2 \\
   - 37092  &  281 &   3 &  11  &   2 &  4 \\
   -  9348  &   41 &  19 &   3  &   3 &  3 \\
   -255972  &  257 &   3 &  83  &   4 &  2 \\
   -101796  &   17 & 499 &   3  &   2 &  5 \\
            &      &     &      &     &
   \end{array} \quad    \begin{array}{r|rrr|rr}
        d   &  p   &   q &   q' &   m &  n \\ \hline
   - 25764  &  113 &  19 &   3  &   3 &  4 \\
   -132612  &  257 &  43 &   3  &   4 &  3 \\
   - 75108  &  569 &   3 &  11  &   5 &  2 \\
   -169796  &   17 &  11 & 227  &   2 &  6 \\
   - 78276  &  593 &   3 &  11  &   3 &  5 \\
   -329988  &  257 &   3 & 107  &   4 &  4 \\
   -106788  &  809 &   3 &  11  &   5 &  3 \\
   -1886244 & 8273 &  19 &   3  &   6 &  2
\end{array} $$

\bigskip

By Pari, we found that for the cases where $m<n-1$, the class groups
agree with those predicted by Corollary~\ref{C2}:
$$ \begin{array}{r|c}
  d \quad   & \Cl(k(\sqrt{-p}\,)) \\ \hline
   -  37092 & (2,4,80) \\
   - 101796 & (2,4,32) \\
   -169796  & (2,4,64)
   \end{array} $$

\bigskip

\section{{\bf Appendix}}
In this appendix,  we  review   some properties of group commutators and of $p$-groups and introduce some notation; see the exposition given, for example, in \cite{BS1}. Assume for the moment that $G$ is an  arbitrary group written multiplicatively. If $x,y\in G$, then define  $[x,y]=x^{-1}y^{-1}xy$, the commutator of $x$ with $y$. More generally, since $[*,*]$ is not associative, we define (inductively on $\ell$)
$[x_1,\dots,x_\ell]=[[x_1,\dots, x_{\ell-1}],x_\ell],$ for all $x_j\in G$ and $\ell >2$. If $A$ and $B$ are nonempty subsets of $G$, then $[A,B]$ is the subgroup of $G$ given as $[A,B]=\la \{[a,b]:a\in A,b\in B\}\ra.$  For example,  the commutator subgroup $G'$ is defined as $[G,G]$; also $G''=(G')'.$ Recall, too, that a group $G$ is metabelian if $G''=1$, i.e., $G'$ is abelian.

Of particular use in our analysis is the lower central series $\{G_\ell\}$ of $G$, which is defined inductively as:
$G_1=G, \quad \text{and}\quad G_{\ell+1}=[G,G_\ell], \quad \text{for all}\;\; \ell\geq 1.$
In particular, observe that $G_2=G'.$  The lower central series is especially useful when the group $G$ is nilpotent,  i.e., when
the lower central series terminates in finitely many steps at the identity subgroup $ 1$. As is well known, all finite $p$-groups are nilpotent for  any prime $p$.

 Recall that  for any $x,y,z\in G$,
\begin{align}\label{G3}
    &x^y=y^{-1}xy=x[x,y]\,,\\
    &[xy,z]=[x,z][x,z,y][y,z]\,,\notag\\
    &[x,yz]=[x,z][x,y][x,y,z]\,,\notag\\
    &[x,y^{-1}]^y=[x^{-1},y]^x=[y,x]=[x,y]^{-1}\,,\notag\\
    &[x,y^{-1},z]^y[y,z^{-1},x]^z[z,x^{-1},y]^x=1\;\;\;\;(\mbox{{\em Witt Identity}}), \notag\\
    & [x,y,z][y,z,x][z,x,y]\equiv 1 \bmod G''\;\;\;\;(\mbox{{\em Witt congruence}}). \notag
    \end{align}
(cf.\,\cite{Bla},  for example).
\bigskip

We also introduce here some notation:  Suppose $G=\la a_1,a_2,\dots,a_r\ra$; then we will let
$$c_{ij}=[a_i,a_j]\,,\quad c_{ijk}=[a_i,a_j,a_k]\,, \;\;\;\text{etc.}$$

\medskip
As usual, we  let  $G^p=\la x^p:x\in G\ra$.

We now collect some consequences of these properties in the situation needed in this work:
\begin{prop}\label{P3} Let $G=\la a_1,a_2,\dots,a_r\ra$ be a metabelian $p$-group; then
\begin{align*}
    &[a_ia_j,a_\ell]=c_{i\ell}c_{j\ell}c_{i\ell j}\,,\\
    &[a_i,a_ja_\ell]=c_{ij}c_{_i\ell}c_{ij\ell}\,,\\
    & c_{ij\ell}=c_{ji\ell}^{-1}\,,\\
    & c_{ij\ell}c_{\ell ij}c_{ji\ell}=1\,.
    \end{align*}
\end{prop}
The proof follows readily from (\ref{G3}).

\bigskip

Now we recall some facts about the Frattini subgroup of a $p$-group and the Burnside Basis Theorem. (Refer to Kapitel III in \cite{Hup}, for example.)   Let $G$ be a finite $p$-group.
Then recall that the Frattini subgroup of $G$, $\Phi(G)$, is the subgroup $G^p G'$ of $G$. When $p=2$, as in our case, then $\Phi(G)=G^2$, since $G'\subseteq G^2$.  Here is a version of the Burnside Basis Theorem:

\medskip

\noindent   {\em Let $G$ be a finite $p$-group. Then the elements $x_1,\dots, x_d$ generate $G$ if and only if the cosets $x_1\Phi(G),\dots,x_d\Phi(G)$ generate $G/\Phi(G).$}

\medskip
Recall, too, that the rank $r$ of a finite abelian $p$-group is the
number of nontrivial cyclic summands in its direct sum
decomposition. (More generally, its $p^n$-rank is the number of cyclic
summands of order at least $p^n$.) As a consequence of the Burnside
Basis Theorem we see that for a finite $p$-group $G$, the order
$|G/\Phi(G)|=p^r$ where $r$ is the rank of $G^{\ab}=G/G'$. Again by
the Burnside Basis Theorem, a minimal set of generators of $G$ has
cardinality equal to the rank of the abelianization
$G^{\ab}$. Finally, if $\{x_1G',\dots,x_r G'\}$ is a basis of
$G^{\ab}$, then $\{x_1,\dots,x_r\}$ is a minimal generating set of
$G$.

\bigskip
\noindent {\em Clarification:} With respect to \cite{BLS2}, p.\,25, we
now give the following less ambiguous definition of {\em primitive}: A
{\em primitive} solution of $d_1d_2X^2+d_3Y^2=Z^2$ is a nontrivial
rational solution $(x,y,z)\in\Q^3$ such that $\alpha=z+x\sqrt{d_1d_2}$
is integral and not divisible by any rational integers ($\not=\pm 1$).
(Notice that this allows for possible semi-integral rational
solutions.)

\newpage

{\bf Addresses}

\medskip

\author{Elliot Benjamin}

\address{Harold Abel School of Psychology,  Capella University,

  Minneapolis, MN 55402, USA; www.capella.edu}

\email{ben496@prexar.com}

\medskip

\author{Franz Lemmermeyer}

\address{M\"orikeweg 1, 73489 Jagstzell, Germany}

\email{hb3@uni-heidelberg.de}

\medskip

\author{C.\,Snyder}

\address{Department of Mathematics and Statistics, UMaine, Orono, ME 04469, USA}

\email{wsnyder@maine.edu}

\end{document}